\documentclass[final,5p,times,twocolumn,authoryear]{elsarticle}

\usepackage{epsfig}

\usepackage{amssymb}
\usepackage{amsmath}
\usepackage{amsfonts} 
\usepackage{subcaption}

\usepackage{amsthm}
\newtheorem{theorem}{Theorem}[section]
\newtheorem{remark}[theorem]{Remark}

\journal{Journal of Process Control}
\begin{document}

\begin{frontmatter}


\cortext[cor]{Corresponding author}

\title{Comparing Model-based Control Strategies for a Quadruple Tank System: \\
Decentralized PID, LMPC, and NMPC}


\author[COMP]{Anders H. D. Christensen}
\address[COMP]{Department of Applied Mathematics and Computer Science, Technical University of Denmark, Matematiktorvet, Building 303B, Kgs. Lyngby, DK-2800, Denmark}
\author[COMP]{Tobias K. S. Ritschel} 
\author[COMP]{Jan Lorenz Svensen}
\author[COMP,2C]{\\Steen Hørsholt}
\author[KT]{Jakob Kjøbsted Huusom}
\address[KT]{Department of Chemical and Biochemical Engineering, Technical University of Denmark, Søltofts Plads 228A, Kgs. Lyngby, DK-2800, Denmark}
\author[COMP,2C]{John Bagterp Jørgensen\corref{cor}}
\address[2C]{2-control ApS, DK-7400 Herning, Denmark}

\ead{jbjo@dtu.dk}

\begin{abstract}
This paper compares the performance of a decentralized proportional-integral-derivative (PID) controller, a linear model predictive controller (LMPC), and a nonlinear model predictive controller (NMPC) applied to a quadruple tank system (QTS). We present experimental data from a physical setup of the QTS as well as simulation results. The QTS is modeled as a stochastic nonlinear continuous-discrete-time system, with parameters estimated using a maximum-likelihood prediction-error-method (ML-PEM). The NMPC applies the stochastic nonlinear continuous-discrete-time model, while the LMPC uses a linearized version of the same model. We tune the decentralized PID controller using the simple internal model control (SIMC) rules. The SIMC rules require transfer functions of the process, and we obtain these from the linearized model. We compare the controller performances based on systematic tests using both the physical setup and the simulated QTS. We measure the performance in terms of tracking errors and rate of movement in the manipulated variables. The LMPC and the NMPC perform better than the decentralized PID control system for tracking pre-announced time-varying setpoints. For disturbance rejection, the MPCs perform only slightly better than the decentralized PID controller. The primary advantage of the MPCs is their ability to use the information of future setpoints. We demonstrate this by providing simulation results of the MPCs with and without such information. Finally, the NMPC achieves slightly improved tracking errors compared to the LMPC but at the expense of having a higher input rate of movement.
\end{abstract}



\begin{keyword}
Quadruple Tank System 
\sep Decentralized PID
\sep Linear MPC
\sep Nonlinear MPC
\sep Parameter estimation
\sep Experiments and simulations
\sep Anticipatory control
\end{keyword}



\end{frontmatter}


\section{Introduction}
\label{sec:introduction}

\begin{figure}[tb]
    \centering
    \includegraphics[width = 0.48\textwidth]{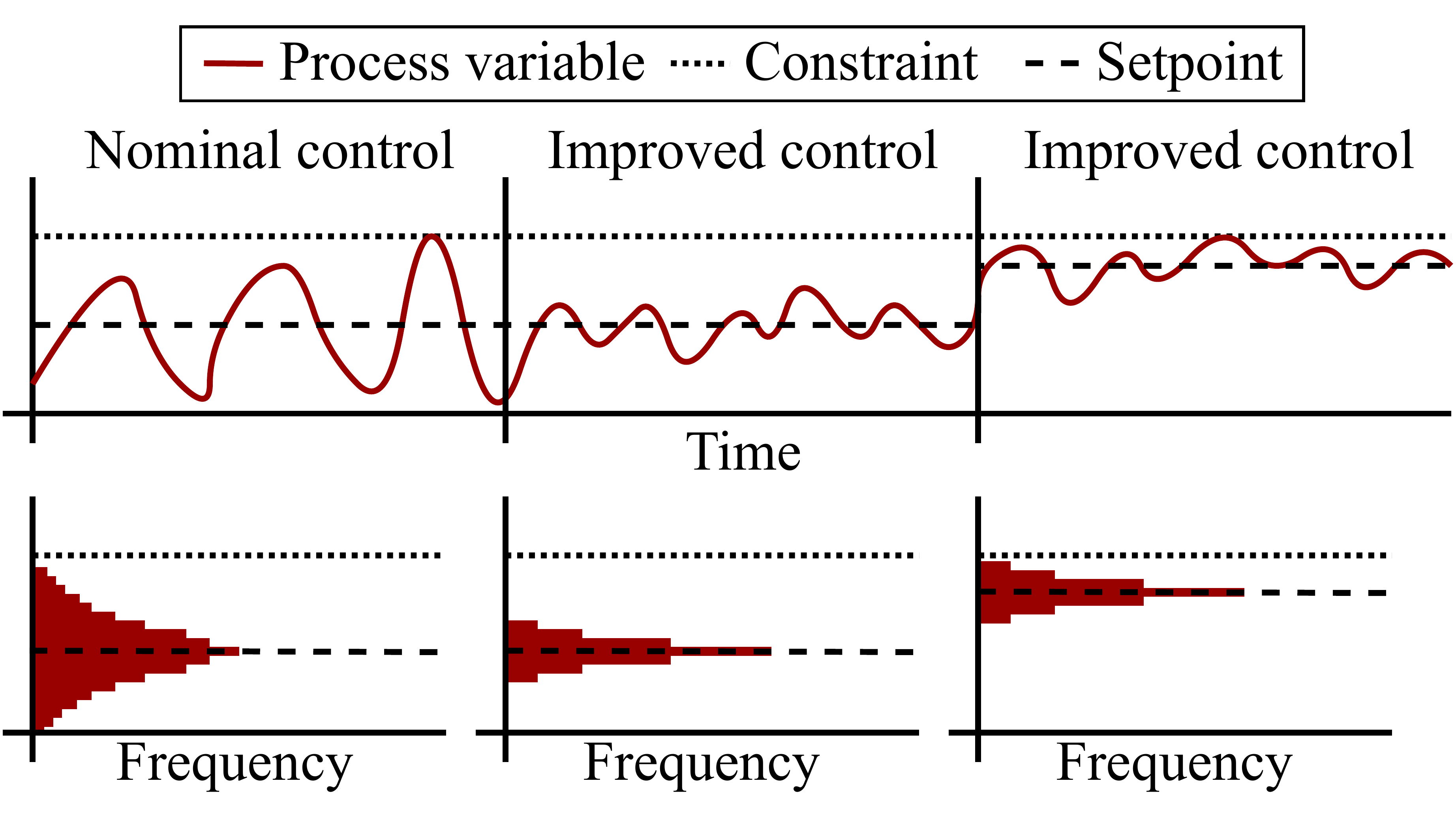}
    \caption{Squeeze and shift principle: Improving the nominal control system squeezes the process variability, allowing operation to be shifted closer to the constraint.}
    \label{fig:squeeze_and_shift}
\end{figure}

In many industrial processes, optimal operating conditions are close to system constraints. Advanced control strategies are often used to operate industrial processes because they can take constraints into account and reduce the effects of disturbances compared to simpler control systems. Figure \ref{fig:squeeze_and_shift} illustrates this situation. Here,  advanced control techniques improve a nominal control system by squeezing the process variability and shifting the average operation closer to the constraint \citep{Seborg:etal:2011,Stephanopoulos:1984,Stephanopoulos:2025,Ogunnaike:Ray:1994,Marlin:2000,Bequette:2003,Brosilow:Joseph:2002}.

The two main approaches for advanced control in the process industries are model predictive control (MPC) and advanced regulatory control (ARC) techniques \citep{Skogestad:2023}. Collectively, MPC and ARC may be referred to as advanced process control (APC). MPC is widely used in the process industries due to its ability to handle constraints in complex multivariable systems \citep{Qin:Badgwell:2000,Qin:Badgwell:2003, Bauer:Craig:2008, Maciejowski:2002, Rawlings:Mayne:Diehl:2022, Camacho:Bordons:1995,Borelli:Bemporad:Morari:2017}. An MPC algorithm is a feedback strategy. It solves a sequence of open-loop optimal control problems (OCPs) in a moving horizon implementation. On the other hand, an ARC system is a decomposed control system that uses proportional-integral-derivative (PID) strategies and other control elements. Decomposing control systems into simple ARC elements provides a powerful technique for controlling complex multivariable systems with constraints. ARC and PID-type controllers are the most commonly used feedback control systems in the process industries \citep{Skogestad:2023, Astroem:Murray:2009, Astrom:Hagglund:1995, Visioli:2006}. ARC and PID-type controllers are also used for feedforward control in the process industries \citep{Guzman:Hagglund:2024}.

\begin{figure}[tb]
    \centering
    \includegraphics[width = 0.48\textwidth]{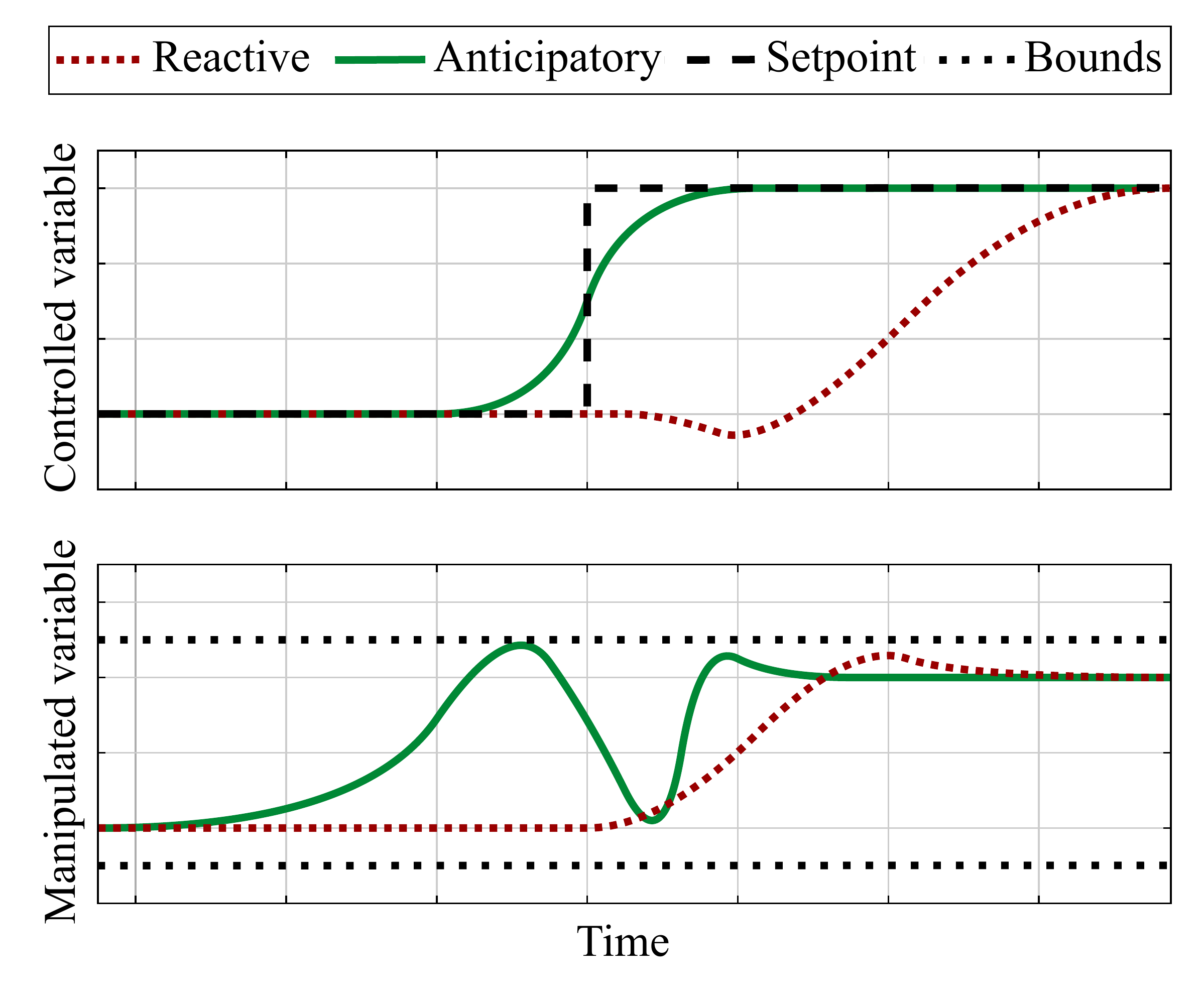}
    \caption{Illustration of the reactive and anticipatory behavior in control systems.}
    \label{fig:anticpatory_action}
\end{figure}

The anticipatory behavior of the MPC is undoubtedly one of the main advantages compared to PID-type control systems. As illustrated in Figure 
\ref{fig:anticpatory_action}, the anticipatory behavior offers better tracking capabilities for time-varying setpoints than reactive-based control systems, when future setpoint changes are known in advance. Nonlinear MPC (NMPC) can further improve setpoint tracking compared to linear MPC (LMPC) for processes with strong nonlinear dynamics \citep{Kamel:etal:2017}. This is because nonlinear control handles process nonlinearities in a wider range of operations compared to linear control strategies \citep{Slotine:Li:1991}. However, NMPC is still less widely used in the process industries compared to LMPC \citep{Qin:Badgwell:2003,Bauer:Craig:2008,Forbes:etal:2015,Darby:2021,Badgwell:Qin:2021,Niu:Xiao:2022}. The polymer manufacturing industry is an exception. In this industry, NMPC has been applied extensively for grade transitions \citep{Henson:1998,Qin:Badgwell:2003,Naidoo:etal:2007,Bindlish:2015,Skaalen:etal:2016}.

Due to the use of a process model in an MPC algorithm, any plant-model mismatch will affect control performance and may result in offsets from setpoints. PID-type control systems naturally handle this offset due to their integral action. To achieve a similar feature in an MPC algorithm, a state estimator with integrating disturbance models used to estimate the unknown disturbances is usually included \citep{Muske:Badgwell:2002, Pannocchia:Rawlings:2003, Pannocchia:2015, Maeder:Morari:2010, Morari:Maeder:2012,Turan:etal:2024,Kuntz:Rawlings:2024}. In these models, the disturbance process is assumed to be integrated white noise. This leads to a deadbeat observer for the disturbance process. Deadbeat observers have erratic behavior and large variance. To obtain offset free observers with better performance than deadbeat observers, the disturbance process may be modeled as a filter with an integrator \citep{Hagdrup:etal:2016,Huusom:etal:2010,Huusom:etal:2011acc,Huusom:etal:2011ifacwc,Huusom:etal:2012,Olesen:etal:2013siso,Olesen:etal:2013mimo}. Generalized predictive control (GPC) also achieves offset-free control by using a filter with an integrator~\citep{Clarke:etal:1987gpc1,Clarke:etal:1987gpc2,Clarke:Mohtadi:1989,Bitmead:Gevers:Wertz:1990}. Offset-free control in MPC may also be achieved using integral control outside the MPC loop \citep{Waschl:2014, Waschl:2015} and by penalizing the integrated error in the optimal control problem \citep{Hendricks:etal:2008,Franklin:etal:1997}. In \cite{Waschl:2014,Waschl:2015}, the offset-free MPC algorithm is decomposed into a nominal disturbance-free MPC loop and an external setpoint adaptation based on proportional-integral (PI) controllers that handle the influence of unknown disturbances. The method described by \citep{Hendricks:etal:2008} and \cite{Franklin:etal:1997} also uses a nominal disturbance-free model in the MPC. The key difference between the disturbance-model approach and these other approaches for offset free control is whether a deliberate model error (the disturbance model) is included in the filtering and prediction model. In the present paper, we model the unknown disturbances as continuous-time Brownian motion (integrated white noise) and the dynamical models describing the plant are systems of stochastic differential equations.  

An ARC system consisting of several PID-type controllers requires many tuning parameters. Finding PID tuning parameters that give similar performance to an MPC can be a challenging task \citep{Pannocchia:etal:2005}. The internal model control (IMC) design procedure offers simple model-based tuning methods for PID-type control systems with only the closed-loop time constants as tuning parameters \citep{Rivera:Morari:Skogestad:1986, Morari:Zafiriou:1989, Skogestad:2003}. 
As a result, the closed-loop performance between a well-tuned ARC system and an MPC might not be so different for situations where it is impossible to supply information regarding future events. It is even shown by \cite{Soroush:Muske:2000} that some MPC algorithms with a prediction horizon of one are similar to a PID-type controller with optimal integrator wind-up. However, for situations where future setpoints can be provided, MPC strategies are usually better than ARC systems \citep{Skogestad:2023}.

Test systems are frequently used to demonstrate multivariable control concepts and compare different control strategies' performance. The quadruple tank system (QTS) \citep{Johansson:2000} is a widely used test system in the research community. \cite{Johansson:2000} applies a manually tuned decentralized PI control system to control the water levels in the two bottom tanks for a physical setup of the QTS. However, more advanced PID-type control systems with auto-tuning capabilities are also applied to the QTS \citep{Ionescu:etal:2016}.
Some research papers present control strategies that address model uncertainties as part of the control design \citep{Huusom:etal:2007, Gurjar:etal:2021, Shah:Patel:2019}. \cite{Huusom:etal:2007} present inventory controllers using a data-driven iterative feedback strategy to control the QTS subject to model-plant mismatch. \cite{Gurjar:etal:2021} use robust control in the form of sliding mode control techniques to handle unknown valve constants, while \cite{Shah:Patel:2019} use sliding mode control to handle time delays.  \cite{Yu:etal:2016} compare the performance of constrained $\mathcal{H}_{\infty}$-control based on linear matrix inequalities (LMIs) for the QTS to unconstrained $\mathcal{H}_{\infty}$-control and decentralized PID-controllers. \cite{Doyle:etal:1999,Doyle:etal:2000} and \cite{Gatzke:etal:2000}  provide simulation and experimental test of internal model control (IMC) based on co-prime factorizations that are related to the Youla-Kucera parametrization \citep{Mahtout:etal:2020,Kucera:2011,Tay:etal:1998,Morari:Zafiriou:1989}. There are also several research papers demonstrating MPC-based strategies for the QTS. \cite{Sazuan:Joergensen:2018} present unconstrained and constrained LMPCs as well as the achievable performance for different QTS designs. \cite{Nirmala:etal:2011} and \cite{Kumar:etal:2019} discuss LMPC for the QTS, while \cite{Zong:etal:2010} compare fast NMPC and constrained LMPC algorithms for the QTS. \cite{Hedengren:etal:2014cace} present NMPC including system identification for a simulated QTS. \cite{Raff:etal:2006} consider the effect of NMPC with stability constraints for a simulated as well as an experimental QTS. \cite{Santos:etal:2012} propose a robust discrete-time LMPC algorithm with explicit dead-time compensation. \cite{Mercangoz:Doyle:2007}, \cite{Alvarado:etal:2011}, \cite{Segovia:etal:2019,Segovia:etal:2021}, and \cite{Grancharova:Johansen:2018} present distributed LMPC and NMPC algorithms and apply these to the QTS. In addition, \cite{Blaud:etal:2022} compare LMPC, linear time-varying MPC, and NMPC with a neural-network-based MPC strategy using the QTS as the test system. They do this to demonstrate the influence of different prediction models in MPC strategies.

Research papers comparing PID-type control systems with MPCs using simulations of other test systems include \cite{Huang:Riggs:2002}, \cite{Taysom:etal:2017}, and \cite{Petersen:etal:2017}. \cite{Huang:Riggs:2002} present a comparative study of a PI controller with an MPC for controlling a gas recovery unit. They show that the MPC provides significant economic advantages over the PI controller.  \cite{Taysom:etal:2017} compare two MPCs with two well-tuned PID-type control systems for temperature control in friction stir welding. They recommend using PID-type control for situations where the temperature is maintained at the setpoint and disturbances are unknown or not modeled. However, they suggest the MPCs when disturbances can be anticipated. \cite{Petersen:etal:2017} compare PI control against an MPC with real-time optimization and economic NMPC for the optimization of spray dryer operation. It is shown that both MPC algorithms increase the profitability of the operation. 
However, to our knowledge the open literature contain no systematic comparisons between PID-type control systems and MPCs for a physical setup of the QTS.  

In this paper, we demonstrate the closed-loop performance differences between a well-tuned PID-type control system, an LMPC, and an NMPC. We do this by comparing these three control systems' capabilities to track predefined time-varying setpoints for a physical setup of the QTS. We use a decentralized approach for the PID-type control system, i.e., we implement two independent single-input-single-output (SISO) PID control loops. Each SISO PID control loop pairs one manipulated variable (MV) with one controlled variable (CV). This decentralized PID control system applies the simple internal model control (SIMC) tuning rules \citep{Skogestad:2003,Skogestad:Postlethwaite:2005}. SIMC are model-based tuning rules. We implement integrator windup to handle constraints in the MVs. The LMPC combines a linear input-bound constrained OCP and a continuous-discrete Kalman filter (CD-KF) for estimating the states and unmeasured disturbances. For the NMPC, we combine a nonlinear OCP with input-bound constraints and a continuous-discrete extended Kalman filter (CD-EKF) for the estimation of states and unmeasured disturbances. We model the QTS as a nonlinear stochastic continuous-discrete-time system. We reduce the plant-model mismatch by estimating the parameters in the model using a maximum-likelihood prediction-error-method (ML-PEM). We use this model directly in the NMPC, and we apply a linearized version of that model in the LMPC. Transfer function representations of the linearized model is used in the SIMC tuning rules for the decentralized PID control system.

We compare the performance of the controllers in terms of tracking errors and rate of movement in the MVs. We do this using both experimental and simulated data. The real-time control software used for the experimental studies is constructed without use of commercial software as described by \cite{Andersen:etal:ecc:2023}. The real-time software for implementation of high-level controllers (LMPC, NMPC, decentralized PID) is based on fundamental programming principles for timers, network communication, and scientific computing for systems and control \citep{Burns:Wellings:2009,Williams:2006,Abramovitch:etal:2023accRTStutorial,Abramovitch:2015cca,Fiedler:etal:2023CEPdompc,Lucia:etal:2017CEPdompc,Christensen:Jorgensen:2024}. By the simulation studies, we demonstrate that a key advantage of the MPC algorithms is their ability to utilize future setpoint information. The simulation studies emulate the real-world scenarios for MPC with and without future setpoint information. The simulation results show, that the MPCs with future setpoint information perform better than the decentralized PID control system. The decentralized PID control system performs slightly better than the MPCs in terms of tracking errors when only current setpoint information is available.

This paper extends the work in \cite{Andersen:etal:cpc:2023} by explicitly showing the three control algorithms, including a detailed description of disturbance modeling in the state estimators. We provide a detailed explanation of how the parameter estimation scheme estimates the model and the Kalman filter parameters. Finally, we provide additional simulation studies that demonstrate the performance difference between the three controllers in various scenarios.

The remaining part of the paper is organized as follows. 
Section \ref{sec:modeling} presents nonlinear and linearized models for the QTS. Section \ref{sec:state_estimation} introduces the CD-KF and CD-EKF including disturbance-augmentation for offset-free estimation and control. Section \ref{sec:parameter_estimation} presents the parameter estimation scheme used for estimating the parameters in the model of the QTS. In Section \ref{sec:control}, we describe the procedure for design of the decentralized PID control system, the LMPC, and the NMPC. Section \ref{sec:results_discussion} presents and discusses the experimental and simulated results.  Finally, we present conclusions in Section \ref{sec:conclusions}.

\section{Modeling the quadruple tank system}
\label{sec:modeling}

\begin{figure}[tb]
    \centering
    \captionsetup{justification=centering}
    
    \begin{subfigure}[b]{0.13\textwidth}
        \centering
        \includegraphics[width=\textwidth]{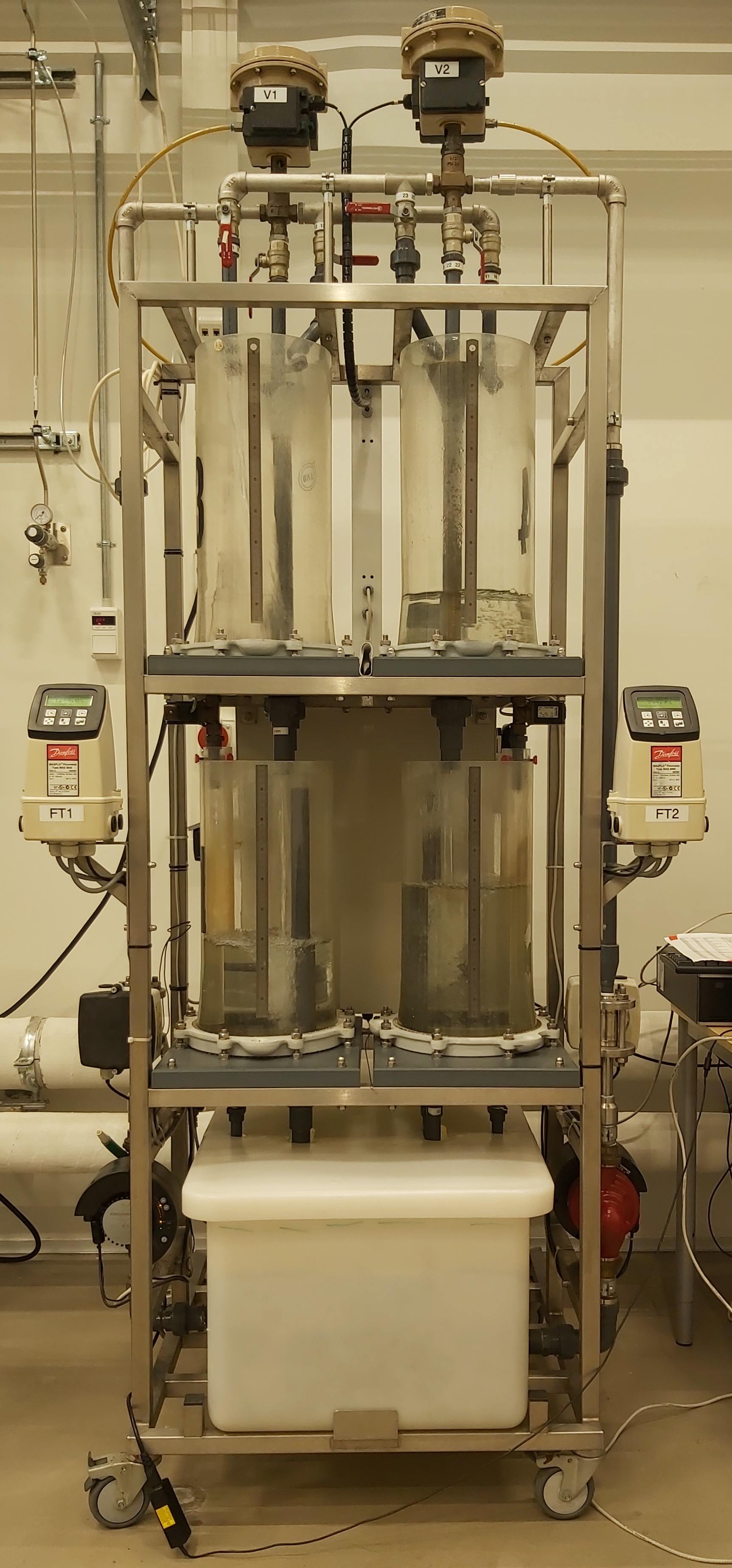}
        \caption{Physical setup.}
        \label{fig:QTS_real}
    \end{subfigure}
    \hfill
    \begin{subfigure}[b]{0.3\textwidth}
        \centering
        \includegraphics[width=\textwidth]{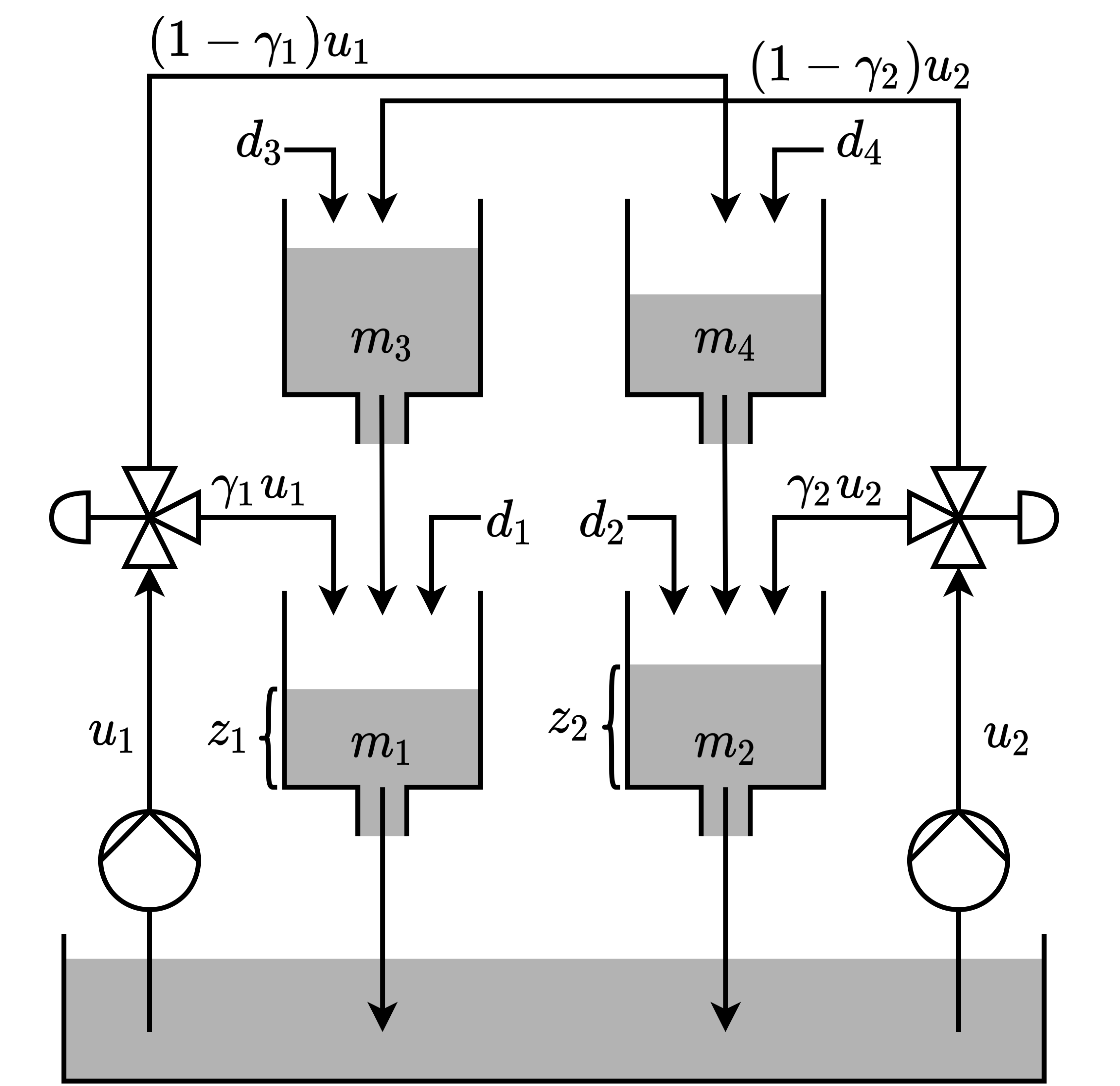}
        \caption{Schematic diagram \citep{Andersen:etal:cpc:2023}.}
        \label{fig:QTS_diagram}
    \end{subfigure}
    \caption{The quadruple tank system.}
    \label{fig:QTS}
\end{figure}

We consider the QTS shown in Figure \ref{fig:QTS}.
The QTS consists of four interconnected water tanks. The MVs are the two pumps that fill these tanks with water from a reservoir. Pump 1 fills tank 1 and tank 4 and pump 2 fills tank 2 and tank 3. Two valves control the distribution of water flow to these tank, i.e., valve 1 controls the fraction of flow from pump 1 to tank 1, and valve 2 controls the flow from pump 2 to tank 2. We measure the water levels in all tanks using pressure sensors. The CVs of the QTS are the water levels in the two bottom tanks. We develop three models for the QTS used in the design of the decentralized PID control system, the LMPC, and the NMPC. 

\subsection{Nonlinear stochastic continuous-discrete-time system}
We model the QTS as a nonlinear stochastic continuous-discrete-time system of the form
\begin{subequations}
\label{eq:model}
\begin{alignat}{1}
        dx(t) &= f(x(t), u(t), d(t), \theta)dt+\sigma( \theta)d\omega(t),\label{eq:SDE}\\
        y(t_k) &= g(x(t_k),\theta)+v(t_k),\label{eq:meas_eq}\\
        z(t) &= h(x(t), \theta).\label{eq:CVs}
\end{alignat}
\end{subequations}
$t$ is time and $x(t) = \begin{bmatrix}
m_1(t); & m_2(t); & m_3(t); & m_4(t)
\end{bmatrix}$ is the state vector representing the masses [g] of water in the tanks. The MVs, $u(t) = \begin{bmatrix}
u_1(t); & u_2(t)
\end{bmatrix}$, are the inflows from
the two pumps to the tanks with units cm$^3$/s. $y(t_k) = \begin{bmatrix}
y_1(t_k); & y_2(t_k); & y_3(t_k); & y_4(t_k)
\end{bmatrix}$ and $z(t) = \begin{bmatrix}
z_1(t); & z_2(t)
\end{bmatrix}$ represent the measured water levels in cm in the tanks and the CVs, respectively. We denote $d(t) = \begin{bmatrix}
d_1(t); & d_2(t); & d_3(t); & d_4(t)
\end{bmatrix}$ as the vector of disturbances. These disturbances represent plant-model mismatch in the form of unknown inflows [cm$^3$/s] in all the tanks. 
$\omega(t)$ is a standard Wiener process, i.e., $d\omega(t) \sim N_{iid}(0, Idt)$ [$\sqrt{\mathrm{s}}$] and $\sigma(\theta)$ represent the state-independent diffusion coefficients [$\text{g}/\sqrt{s}$]. The measurement noise, $v(t_k) \sim N_{iid}(0, R)$, is a discrete-time normally distributed stochastic process with covariance $R~=~\mathrm{diag}\big(\begin{bmatrix}
r_{1}^2 & r_{2}^2 & r_{3}^2 & r_{4}^2
\end{bmatrix}\big)$. Finally, $\theta$ is the parameter vector which we assume is time-invariant.

We model the system of stochastic differential equations in \eqref{eq:SDE} as the mass balances
\begin{alignat}{1}
        d m_i(t) &= \left( \rho q_{i,in}(t) - \rho q_{i,out}(t) \right) dt + \sigma_{i} d\omega_i(t), 
        \label{eq:differential_quadruple_tank_system}
\end{alignat}
for $i \in \{1, 2, 3, 4\}$ with 
$\rho$ = 1.0 g/cm$^3$ being the density of water. The inflows $q_{i,in}(t)$ are described by
\begin{subequations}
    \label{eq:QTS:Inflows}
    \begin{alignat}{1}
    q_{1,in}(t) &= \gamma_1u_1(t)+q_{3,out}(t)+d_1(t),\\
            q_{2,in}(t) &= \gamma_2u_2(t)+q_{4,out}(t)+d_2(t),\\
            q_{3,in}(t) &= (1-\gamma_2)u_2(t)+d_3(t),\\
            q_{4,in}(t) &= (1-\gamma_1)u_1(t)+d_4(t),
    \end{alignat}
\end{subequations}
where $\gamma_1, \gamma_2 \in (0,1)$ represent the fixed valve configurations. $q_{i,out}(t)$ represent the outflows and we describe these as
\begin{subequations}
    \label{eq:QTS:Outflows}
    \begin{alignat}{1}
    q_{i,out}(t) &= a_i\sqrt{2g_ah_i(t)}, \\ h_i(t) &= \frac{m_i(t)}{\rho A_i},\label{eq:water_levels}
    \end{alignat}
\end{subequations}
for $i \in \{1, 2, 3, 4\}$. $g_a = 981\:\mathrm{cm}/\mathrm{s}^2$ is the acceleration of gravity, and $a_i$ and $A_i$ represent the cross-sectional areas of the outlet tubes and tanks, respectively \citep{Johansson:2000}. The CVs are the water levels $h_1(t)$ and $h_2(t)$ from \eqref{eq:water_levels}. Hence, we formulate \eqref{eq:CVs} as
\begin{equation}
   h(x(t), \theta) = C_z(\theta)x(t), \quad C_z(\theta) =\begin{bmatrix}
 \frac{1}{\rho A_1}&0  & 0 & 0\\ 
0 &\frac{1}{\rho A_2}  &  0& 0
\end{bmatrix}.
\end{equation}
Similarly, we express the measurement function in \eqref{eq:meas_eq} as
\begin{equation}
    g(x(t_k), \theta) = C(\theta)x(t_k), 
\end{equation}
where
\begin{equation}
\label{eq:C_matrix}
    C(\theta) = \mathrm{diag}\bigg(
    \begin{bmatrix}
\frac{1}{\rho A_1} & \frac{1}{\rho A_2} & \frac{1}{\rho A_3} & \frac{1}{\rho A_4}
\end{bmatrix}\bigg).
\end{equation}

For simplicity, we model the diffusion coefficient in \eqref{eq:SDE} using the diagonal matrix,
\begin{equation}
    \sigma(\theta) = \mathrm{diag}\big(\begin{bmatrix}
\sigma_{1} & \sigma_{2} & \sigma_{3} & \sigma_{4}
\end{bmatrix}\big).
\end{equation}

\subsubsection{Nominal parameters}

We obtain the nominal parameters of the QTS by measuring the cross-sectional areas as $a_i = 1.13\; \mathrm{cm}^2$ and $A_i = 380.13\; \mathrm{cm}^2$ for $i \in \{1, 2, 3, 4\}$. 
For practical reasons, we choose the valve configurations as $\gamma_j = 0.35$ for $j \in \{1, 2\}$. This ensures that most of the water flows through tank 3 and tank 4, thus avoiding emptying the tanks. As a result, the QTS has non-minimum phase characteristics. 
We combine these process parameters and the diffusion coefficients in the parameter vector $\theta$. However, we do not have information regarding the nominal values of the diffusion coefficients.

\subsection{Linearized models}
We derive a linear model of the QTS by linearizing \eqref{eq:model} at the operating point $(x_s, u_s, d_s, y_s, z_s)$,
\begin{subequations}
\label{eq:linear_state_space_model}
\begin{alignat}{1}
       dX(t) &= \bigg(A(\theta)X(t)+B(\theta)U(t)+E(\theta)D(t)\bigg)dt\nonumber \\
       &+\sigma(\theta)d\omega(t), \label{subeq:sde_linear}\\
        Y_k &= C(\theta)X_k+v_k,\\
        Z(t) &= C_z(\theta)X(t).
\end{alignat}
\end{subequations}
 $X(t) = x(t)-x_{s}, U(t)=u(t)-u_{s}, D(t)=d(t)-d_{s}, Y_k=Y(t_k)=y(t_k)-y_{s}=y_k-y_{s}$, and $Z(t)=z(t)-z_{s}$ represent the deviation of the variables from the steady state operating point $(x_{s},u_{s},d_{s},y_{s},z_{s})$. $A(\theta) = \partial_x f(x_s,u_s,d_s,\theta)$, $B(\theta)=\partial_u f(x_s,u_s,d_s,\theta)$, $E(\theta) = \partial_d f(x_s,u_s,d_s,\theta)$, $C(\theta)=\partial_x g(x_s,\theta)$, and $C_z(\theta)=\partial_x h(x_s,\theta)$ 
are the derivatives of the functions $f$, $g$, and $h$ in \eqref{eq:model}. The drift, $f$, in the QTS model specified by \eqref{eq:differential_quadruple_tank_system}-\eqref{eq:QTS:Inflows} is affine in $u$ and $d$, i.e. it is in the form $f(x(t),u(t),d(t),\theta) = F(x(t),\theta) + B(\theta) u(t) + E(\theta) d(t)$. Consequently, $A(\theta) = \partial_x f(x_s,u_s,d_s,\theta) ) = \partial_x F(x_s,\theta)$ is dependent on the operating point, while $B(\theta)=\partial_u f(x_s,u_s,d_s,\theta)=[\, \gamma_1\,0;\,0\,\gamma_2;\,0 \, 1-\gamma_2;\, 0 \, 1-\gamma_1]$ and $E(\theta) = \partial_d f(x_s,u_s,d_s,\theta) = I$ ($I$ is the identity matrix) are independent of the operating point.
We compute the transfer function matrix of \eqref{eq:linear_state_space_model} from MVs to CVs as
\begin{equation}
    G(s) = \begin{bmatrix}
g_{11}(s) & g_{12}(s) \\ 
g_{21}(s) & g_{22}(s)
\end{bmatrix} = C_z(\theta) (s I - A(\theta))^{-1} B(\theta).
\end{equation}
The cross-coupling in the system is represented using the two second-order transfer functions, $g_{12}(s)$ and $g_{21}(s)$. They are both in the form
\begin{equation}
\label{eq:transfer_functions}
    g(s) = \frac{k}{(\tau_1s+1)(\tau_2s+1)},
\end{equation}
where $\tau_1 \geq \tau_2$ are time constants, and $k$ is the steady-state gain. $g_{11}(s)$ and $g_{22}(s)$ are first-order transfer functions, but we do not use these in this paper.
\section{State and disturbance estimation}
\label{sec:state_estimation}

This section presents the CD-EKF and CD-KF used for state and disturbance estimation in the NMPC and LMPC, respectively. The ML-PEM parameter estimation scheme also requires the CD-EKF \citep{Astrom:1980,Kristensen:etal:2004,Jorgensen:2004, Jorgensen:Jorgensen:PEM:2007,Boiroux:etal:2016,Boiroux:etal:ACC:2019,Boiroux:etal:2019,Sarkka:Solin:2019}.

\subsection{Continuous-discrete extended Kalman filter}

The CD-EKF algorithm consists of a filtering step and a one-step prediction step. At time $t_k$, the filtering step updates the current state estimate based on the measurement, $y_k$, and the measurement equation. The prediction step uses the MVs, $u_k$,  to predict the state at time $t_{k+1}$. 

\textit{Filtering:} The CD-EKF obtains the estimate, $\hat{x}_{k|k}$, and its covariance, $P_{k|k}$, using the information of the previous predicted estimate, $\hat{x}_{k|k-1}$, and covariance, $P_{k|k-1}$, as
\begin{subequations}
\label{eq:filter1}
\begin{alignat}{1}
        \hat{x}_{k|k} &= \hat{x}_{k|k-1}+K_ke_k,\\
        P_{k|k} &= (I-K_kC_k)P_{k|k-1}(I-K_kC_k)'+K_kRK_k',
\end{alignat}
\end{subequations}
where
\begin{subequations}
\label{eq:filter2}
\begin{alignat}{1}
    \hat{y}_{k|k-1} &= g(\hat{x}_{k|k-1}, \theta), \quad     C_k = \frac{\partial g}{\partial x}(\hat{x}_{k|k-1}, \theta), \\
     e_{k} &= y_k-\hat{y}_{k|k-1},  \quad  R_{e,k} = R+C_kP_{k|k-1}C_k',\\
    K_k &= P_{k|k-1}C_k'R_{e,k}^{-1}. 
\end{alignat}
\end{subequations}

\textit{Prediction}: 
Given the estimate, $\hat{x}_{k|k}$, and covariance, $P_{k|k}$, the one-step prediction of the state-covariance pair,
\begin{subequations}
\begin{alignat}{1}
        \hat{x}_{k+1|k} &= \hat{x}_{k}(t_{k+1}),\quad
        P_{k+1|k} = P_{k}(t_{k+1}),
\end{alignat}
\end{subequations}
is obtained by numerical solution of
\begin{subequations}
\label{eq:timeupdate_CDEKF}
\begin{alignat}{3}
        \frac{d}{dt}\hat{x}_{k}(t) &= f( \hat{x}_{k}(t), u_{k}, d_{k}, \theta),
    \label{eq:CDEKF_timeupdate}\\
        \frac{d}{dt}P_{k}(t) &= A_{k}(t)P_{k}(t)+P_{k}(t)A_{k}(t)' +\sigma(\theta) \sigma(\theta)' \label{eq:CDEKF_timeupdate_cov},
\end{alignat}
\end{subequations}
for $t \in [t_{k}, t_{k+1}]$ with the initial conditions
\begin{equation}
    \hat{x}_{k}(t_{k}) = \hat{x}_{k|k},\quad
        P_{k}(t_{k}) = P_{k|k}.
\end{equation}
We compute the Jacobian $A_k(t)$ as
\begin{equation}
\label{eq:Akt}
    A_{k}(t) = \frac{\partial f}{\partial x}(\hat{x}_{k}(t), u_{k}, d_{k}, \theta).
\end{equation}
We solve \eqref{eq:timeupdate_CDEKF} using a fourth-order explicit Runge-Kutta scheme with 10 fixed integration steps in the interval $[t_k, \, t_{k+1}]$. 

\subsection{Disturbance estimation}
Solving \eqref{eq:timeupdate_CDEKF} requires $d_k$. We estimate $d_k$ by augmenting \eqref{eq:model} with a disturbance model. We consider the integrating disturbance model
\begin{equation}
\label{eq:distubance_models}
dd_i(t) = \sigma_{d,i}(\theta)d\omega_j(t),
\end{equation}
for $i \in \{1, 2, 3, 4\}$ and $j \in \{5, 6, 7, 8\}$ \citep{Jorgensen:2007,Joergensen:ACC:2007,Pannocchia:Rawlings:2003, Odelson:Rawlings:2003}. We use this simple disturbance model, as we have no prior information regarding the disturbance dynamics. The CD-EKF estimates the states and disturbances simultaneously, by applying the filtering and prediction scheme to the disturbance-augmented model

\begin{subequations}
\label{eq:dist_augmented_model}
\begin{alignat}{1}
        dx_a(t) &= f_a(x_a(t), u(t), \theta)dt+\sigma_a( \theta)d\omega_a(t),\\
        y(t_k) &= g_a(x_a(t_k),\theta)+v(t_k).
\end{alignat}
\end{subequations}
$x_a(t) = [x(t); d(t)]$ is the combined state and disturbance vector and $\omega_a(t)$ represents the vector of Wiener processes for \eqref{eq:model} and \eqref{eq:distubance_models}. We formulate the drift coefficient, $f_a(\cdot)$, the diffusion coefficient, $\sigma_a(\theta)$, and the measurement equation as
\begin{subequations}
\label{eq:disturbance_augmented_model}
    \begin{alignat}{1}
    f_a(x_a(t), u(t), \theta) &=
    \begin{bmatrix}
 f(x(t), u(t), d(t), \theta)\\ 
0
\end{bmatrix},\\
\sigma_a(\theta) &= \mathrm{diag}([\sigma_1\quad\sigma_2\quad \sigma_3\quad \sigma_4\quad \nonumber \\
&\quad\sigma_{d,1}\quad \sigma_{d,2}\quad\sigma_{d,3}\quad \sigma_{d,4}]),\\
g_a(x_a(t_k),\theta) &= g(x(t_k),\theta).
\end{alignat}
\end{subequations}
From the filtering step, the CD-EKF obtains the combined estimate of the state and disturbance vector
\begin{equation}
    \hat{x}_{a,k|k} = [\hat{x}_{k|k}; \,\hat{d}_{k|k}],
\end{equation}
with the covariance $P_{a,k|k}$.

\subsection{Continuous-discrete Kalman filter}

The CD-KF applies the filtering and prediction steps to the disturbance-augmented model \eqref{eq:dist_augmented_model} linearized at the operating point $(x_s, u_s, d_s)$. The CD-KF obtains the state-covariance pair, $\hat{X}_{a,k+1|k}=\hat{X}_{a,k}(t_{k+1})$ and $\hat{P}_{a,k+1|k}=\hat{P}_{a,k}(t_{k+1})$, by solving 
\begin{subequations}
\label{eq:timeupdate_CDKF}
\begin{alignat}{3}
    \frac{d}{dt}\hat{X}_{a,k}(t) &= A_a(\theta)\hat{X}_{a,k}(t)+B_a(\theta)U(t),
    \label{eq:CDKF_timeupdate}\\
    \begin{split}
    \frac{d}{dt}P_{a,k}(t) &= A_a(\theta)P_{a,k}(t)+P_{a,k}(t)A_a(\theta)' 
    +\sigma_a(\theta) \sigma_a(\theta)',
    \end{split}
    \label{eq:CDKF_timeupdate_cov}
\end{alignat}
\end{subequations}
for $t \in [t_{k}, t_{k+1}]$ with the initial conditions 
\begin{equation}
    \hat{X}_{a,k}(t_k) = \hat{X}_{a,k|k}, \quad \hat{P}_{a,k}(t_k) = \hat{P}_{a,k|k}.
\end{equation}
$X_a(\theta) = [X(t); D(t)]$ is the combined states and disturbances in deviation variables and we construct the matrices $A_a(\theta)$ and $B_a(\theta)$ as
\begin{equation}
    A_a(\theta) = \begin{bmatrix}
A(\theta) & E(\theta)\\ 
0 & 0
\end{bmatrix}, \quad B_a(\theta) = \begin{bmatrix}
B(\theta)\\ 
0
\end{bmatrix}.
\end{equation}

The filtering step of the CD-KF computes $\hat{X}_{a,k|k} = [\hat{X}_{k|k}; \,\hat{D}_{k|k}]$ using \eqref{eq:filter1}-\eqref{eq:filter2} with the innovation $e_{k} = Y_k-C_k\hat{X}_{a,k|k-1}$ and $C_k = C(\theta)$ precomputed.

\section{Parameter estimation}
\label{sec:parameter_estimation}

This section presents the ML-PEM method for estimating the parameters in \eqref{eq:model}. 
We apply this method to an estimation data set from the physical setup of the QTS. We evaluate the estimated parameters by simulating \eqref{eq:model} with nominal and estimated parameters and we compare the simulated data with both the estimation data and a validation data set. Finally, we apply the ML-PEM method to estimate the measurement noise covariance and diffusion coefficients in the disturbance-augmented CD-EKF.

\subsection{A maximum-likelihood prediction-error-method}

We estimate the parameters in \eqref{eq:model} using the ML-PEM approach presented by \cite{Kristensen:etal:2004}.
Given a data set of $N$ measurements and MVs, 
\begin{subequations}
\begin{alignat}{3}
Y_N &= \begin{bmatrix}
    y_1,& y_2, &y_3,  & \cdots, & y_N
    \end{bmatrix},\\
    U_N &= \begin{bmatrix}
    u_1,& u_2, &u_3,  & \cdots, & u_N
    \end{bmatrix},
\end{alignat}
\end{subequations}
the maximum-likelihood estimate of the parameter $\theta$, $\theta^*_{ML}$, is the vector of parameters that maximizes the likelihood function $p(Y_N|\theta)$. We construct the likelihood function as the product of conditional densities in the form
\begin{equation}
\label{eq:likelihood_fnc}
    \begin{split}
        p(Y_N|\theta) &= \prod_{k=1}^Np(y_k|Y_{k-1},\theta) \\
        &=  \prod_{k=1}^N\frac{\exp{\bigg(-\frac{1}{2}(y_k-\mu_k)'\Sigma_k^{-1}(y_k-\mu_k)\bigg)}}{(2\pi)^{n_y/2}\sqrt{\det(\Sigma_k)}},
    \end{split}
\end{equation}
where $n_y$ is the dimension of $y_k$, and $\mu_k$ and $\Sigma_k$ are the mean and covariance, respectively. The CD-EKF provides the mean estimate $\hat y_{k|k-1} = E\{y_k|Y_{k-1},U_{k-1},\theta\}$ and the innovation covariance $R_{e,k} = V\{y_k|Y_{k-1},U_{k-1},\theta\}$. Using the innovation, $e_k = y_k-\hat y_{k|k-1}$, we express \eqref{eq:likelihood_fnc} in terms of the CD-EKF as
\begin{equation}
    p(Y_N|\theta) = \prod_{k=1}^N \frac{\exp\bigg(-\frac{1}{2}e_k' R_{e,k}^{-1}e_k\bigg)}{(2\pi)^{n_y/2}\sqrt{\det( R_{e,k})}}.
\end{equation}
We obtain $\theta_{ML}^*$ by solving
\begin{equation}
\label{eq:ML_NLP}
\theta_{ML}^*=\arg\min V_{ML}(\theta),
\end{equation}
with the objective function $V_{ML}(\theta) =-\ln(p(Y_N|\theta))$,
\begin{equation}
     V_{ML}(\theta)= \frac{1}{2}\sum_{k=1}^{N}\bigg(\ln{\det( R_{e,k})}+e_k'R_{e,k}^{-1}e_k \bigg) + \frac{Nn_y}{2}\ln{2\pi}.
\end{equation}

\begin{remark}[Maximum a posteriori estimation]
We could formulate the PEM using a maximum a posteriori (MAP) estimation scheme, since the nominal parameters of \eqref{eq:model} are available \citep{Kristensen:etal:2004}. However, this approach requires the distribution of these parameters, which is unknown. Alternatively, one could manually choose the distribution of the nominal parameters, thus incorporating a tuning mechanism into the PEM scheme that biases the estimates toward the nominal values.
\end{remark}

\subsection{Estimating the parameters in the QTS}

We estimate the parameters of \eqref{eq:model} in two stages. In the first stage, we estimate the measurement noise covariance. In the second stage, we estimate the parameters in \eqref{eq:SDE} using the estimated measurement noise covariance in the CD-EKF. 

\subsubsection{Obtaining data}
We use the physical setup of the QTS to generate estimation and validation data. We do this by storing process data obtained for step changes in the MVs. 
These step changes are chosen such that the data sets contain both step changes in one MV at a time and simultaneous changes in the MVs. No external disturbances affect the physical setup of the QTS in these data sets.

\subsubsection{Obtaining measurement noise covariance}
\label{subseq:obtain_meas_noise_cov}

Before estimating the model parameters, we obtain an estimate of the measurement noise covariance, $R$, required by the CD-EKF. We do this by measuring the distribution of the water level measurements when the QTS is in steady-state. We do this for multiple steady-states by choosing different combinations for the MVs. The measurement distribution in the two upper tanks varied significantly. It is assumed that the main reason for this is that the pressure sensors do not measure correctly due to the presence of turbulence near them. 
Therefore, we multiplied the variances by 1000. Consequently, the measurements from the lower tanks are preferred for parameter estimation.

\subsubsection{Applying the ML-PEM scheme to QTS data}
\indent

We estimate the cross-sectional areas of the outlet tubes and tanks, the fixed valve configurations, and the diffusion coefficients from the estimation data set.
We use the non-augmented model \eqref{eq:model} as the model in the CD-EKF assuming no disturbances, i.e., $d(t) = 0$. As a result, we obtain $\theta_{ML}^*$ that minimizes the plant-model mismatch. 
Table \ref{tab:all_parameters} presents $\theta_{ML}^*$ together with the nominal parameters. We see a significant difference between nominal and estimated values for the cross-sectional areas of the two upper tanks. As discussed in Section \ref{subseq:obtain_meas_noise_cov}, these results are related to the assumption that the pressure sensors did not measure the correct water levels in the upper tanks. The estimated value of valve 1 is also very different from its nominal value.

\begin{table}[tb]
\centering
\caption{Nominal and estimated parameters for \eqref{eq:model}. The measurement noise covariance is pre-estimated and manually tuned.}
\begin{tabular}{crrc}
\hline
Parameter & Nominal & Estimated & Unit \\ \hline
$a_1$  & 1.13 & 1.01  & cm$^2$  \\
$a_2$  & 1.13 & 1.25  & cm$^2$\\
$a_3$  & 1.13 & 1.32 & cm$^2$ \\
$a_4$  & 1.13 & 1.55 & cm$^2$ \\
$A_1$  & 380.13 & 379.84& cm$^2$ \\
$A_2$  & 380.13 & 378.03& cm$^2$\\
$A_3$  & 380.13 & 466.30& cm$^2$\\
$A_4$  & 380.13 & 523.12& cm$^2$\\
$\gamma_1$  & 0.35 & 0.260 & -- \\
$\gamma_2$  & 0.35 & 0.353 & --\\
$\sigma_{1}$ & - & 10.07 $\cdot 10^{-3}$ & g/$\sqrt{\mathrm{s}}$\\
$\sigma_{2}$   & - & 13.09 $\cdot 10^{-3}$ & g/$\sqrt{\mathrm{s}}$ \\
$\sigma_{3}$   & - & 12.50 $\cdot 10^{-3}$ & g/$\sqrt{\mathrm{s}}$\\
$\sigma_{4}$   & - & 16.62 $\cdot 10^{-3}$ & g/$\sqrt{\mathrm{s}}$\\
\hline
\end{tabular}
\label{tab:all_parameters}
\end{table}

\begin{figure}[tb]
    \centering
    \includegraphics[width = 0.5\textwidth]{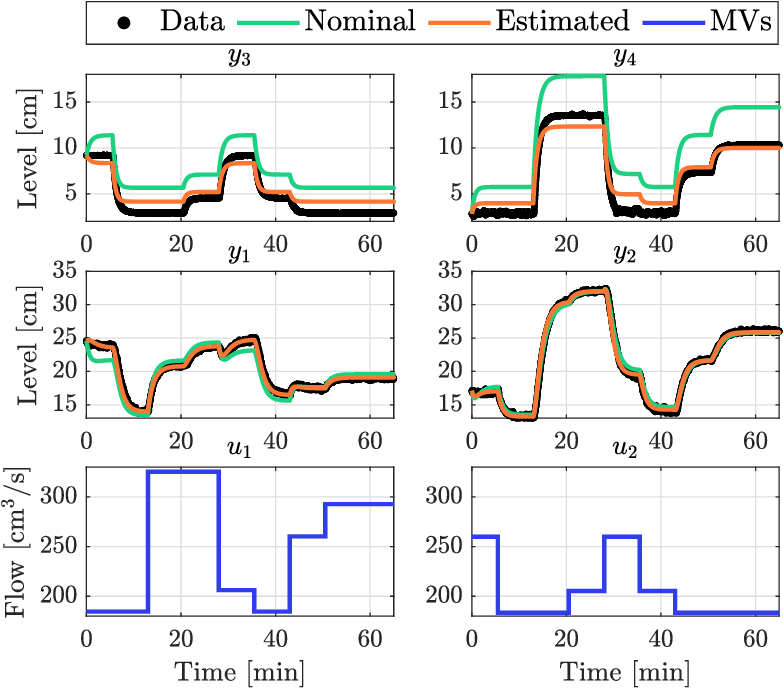}
    \caption{Data used for estimation and simulations with nominal and estimated parameters.}
\label{fig:estimation_data}
\end{figure}

\begin{figure}[tb]
        \centering
    \includegraphics[width = 0.5\textwidth]{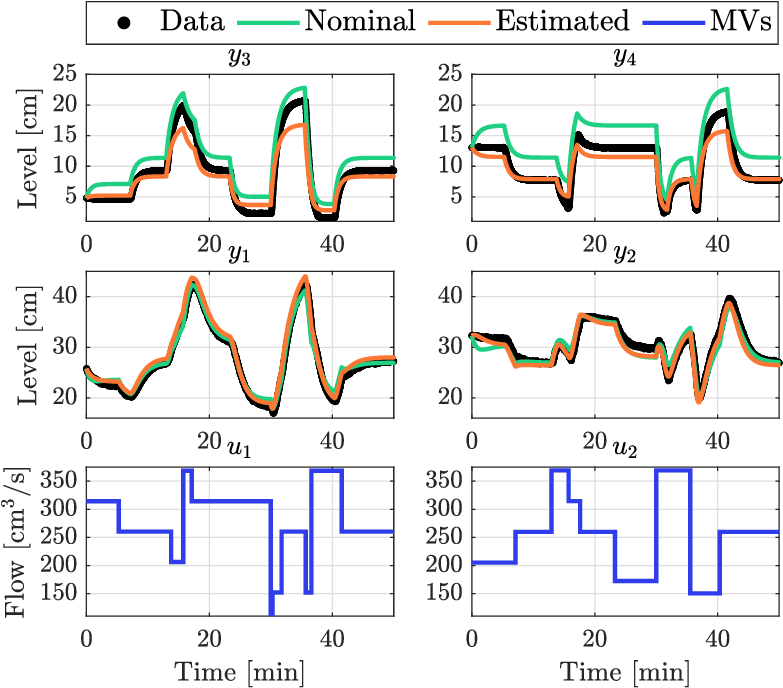}
    \caption{Data used for validation and simulations with nominal and estimated parameters.}
\label{fig:validation_data}
\end{figure}

\subsubsection{Evaluation of estimated parameters}
We evaluate the models by simulating them with estimated and nominal parameters. Figure \ref{fig:estimation_data} shows the estimation data plotted together with open-loop simulations of \eqref{eq:model}. Figure \ref{fig:validation_data} presents simulations compared to the validation data set. 
We measure the goodness-of-fit (GOF) between data from the QTS and the water levels from the open-loop simulations of the system \eqref{eq:model} by computing the averaged normalized root mean squared error,
\begin{equation}
    \mathrm{GOF} = \frac{1}{n_y} \sum_{i=1}^{n_y}\sum_{k=1}^{N} \bigg(1-\frac{||y_{i,k}-\tilde{y}_i(t_k)||}{||y_{i,k}-\mathrm{mean}_k(y_i)||}\bigg)100.
\end{equation}
$y_{i,k}$ and $\tilde{y}_i(t_k)$ for $i \in \{1, 2, 3, 4\}$ are data from the QTS and simulated water levels in tanks using \eqref{eq:SDE}-\eqref{eq:meas_eq} without noise, respectively, and $\mathrm{mean}_k(y_i) = \frac{1}{N}\sum_{k=1}^N(y_{i,k})$.  Table \ref{tab:GOF_estimation_and_validation} presents the GOF for simulations using \eqref{eq:model} with the nominal parameters and the estimates of the parameters. As expected, the GOF is significantly higher when using estimates of the parameters instead of the nominal parameters.

\begin{table}[b]
\centering
\caption{GOF for estimation and validation data.}
\begin{tabular}{llllll}
\hline
Parameters ($\theta$) &  Estimation GOF &  Validation GOF \\ \hline
 Nominal & 47.91\% & 57.28\%  \\
Estimated & 80.41\% & 74.20\%\\
\hline
\end{tabular}
\label{tab:GOF_estimation_and_validation}
\end{table}

\subsection{Tuning of the disturbance-augmented CD-EKF}
The disturbance-augmented CD-EKF requires the measurement noise covariance, $R$, and the diffusion coefficients, $\sigma_a(\theta)$.
We apply the ML-PEM scheme to estimate these in the disturbance-augmented CD-EKF, simultaneously. 
To simplify the procedure, we fix the parameters of the drift term with the estimated parameters in Table \ref{tab:all_parameters} and reuse the estimation data in Figure \ref{fig:estimation_data}. Table \ref{tab:KF_tuning} presents the estimated diffusion coefficients and measurement noise covariance. In this table, the nominal parameters for $R$ are the values obtained through the steady-state procedure described earlier.



\begin{table}[tb]
\centering
\caption{Estimated diffusion coefficients, $\sigma_a(\theta)$, and measurement noise covariance, $R$, using the estimated model parameters of the drift term in Table \ref{tab:all_parameters}.}
\begin{tabular}{crrc}
\hline
Parameter & Nominal & Estimated & Unit \\ \hline
$\sigma_{1}$ & - & 7.25 & g/$\sqrt{\mathrm{s}}$\\
$\sigma_{2}$   & - & 14.92 & g/$\sqrt{\mathrm{s}}$ \\
$\sigma_{3}$   & - & 8.98  & g/$\sqrt{\mathrm{s}}$\\
$\sigma_{4}$   & - & 14.50  & g/$\sqrt{\mathrm{s}}$\\
$\sigma_{d,1}$ & - & 0.47  & g/$\sqrt{\mathrm{s}}$\\
$\sigma_{d,2}$   & - & 3.08  & g/$\sqrt{\mathrm{s}}$ \\
$\sigma_{d,3}$   & - & 3.92  & g/$\sqrt{\mathrm{s}}$\\
$\sigma_{d,4}$   & - & 3.42  & g/$\sqrt{\mathrm{s}}$\\
$r_1^2$          & 0.0349  & 1.44$\cdot 10^{-2}$   & [m$^2$] \\
$r_2^2$          &  0.0340 & 1.34$\cdot 10^{-2}$   & [m$^2$] \\
$r_3^2$          &  0.0025 & 1.00$\cdot 10^{-5}$   & [m$^2$] \\
$r_4^2$          & 0.0065  & 1.00$\cdot 10^{-5}$   & [m$^2$] \\
\hline
\end{tabular}
\label{tab:KF_tuning}
\end{table}

\section{Model-based control algorithm}
\label{sec:control}

This section presents the implementation of the decentralized PID, the LMPC, and the NMPC.
We denote the setpoints for the CVs at time $t_k$ as $\Bar{z}_k~=~\begin{bmatrix}
\Bar{z}_{1,k}; & \Bar{z}_{2,k}
\end{bmatrix}$ and the rate of movement in the MVs as $\Delta u_k = u_{k}-u_{k-1}=\begin{bmatrix}
u_{1,k}-u_{1,k-1}; & u_{2,k}-u_{2,k-1}
\end{bmatrix}$. The sampling time $T_s = t_{k+1}-t_k$ for all three controllers is $T_s=5$ s.

\subsection{Decentralized PID}

The decentralized PID control system consists of two single-input single-output (SISO) PID loops. Due to the valve configurations, pump 1 primarily influences tank 2 and 4 and pump 2 primarily influences tank 1 and 3. Therefore, to ensure good paring between CVs and MVs, we choose the first PID loop to use $y_1$ and $\Bar{z}_1$ to compute $u_2$, and the other PID loop to use $y_2$ and $\Bar{z}_2$ to compute $u_1$. To handle bounds in the MVs, both SISO PID loops include an integrator windup mechanism. Let $(y,z,u)$ be $(y_1,z_1,u_2)$ for PID control loop 1 and $(y_2,z_2,u_1)$ for PID control loop 2. Using this notation, we implement the PID loops as 
\begin{equation}
\label{eq:PID1}
    v_k = \bar{u} + P_k + I_k + D_k, \quad 
         u_k =  
\begin{cases}
    u_{\min} & \text{if } v_k \leq u_{\min}\\
    u_{\max} & \text{if } v_k \geq u_{\max}\\
    v_k,& \text{else }\\
\end{cases},
\end{equation}
where $P_k$, $D_k$, and $I_{k+1}$ are obtained as
\begin{subequations}
\label{eq:PID2}
\begin{alignat}{3}
         &e_k = \Bar{z}_k-y_k,\\
         & s_k = u_k-v_k,\\
       &P_k = K_pe_k,\\
        &D_k = \frac{\tau_d}{\tau_d+NT_s}D_{k-1}-\frac{K_p\tau_dN}{\tau_d+NT_s}\bigg(y_k-y_{k-1}\bigg),\\
        &I_{k+1} = I_k+T_s\frac{K_p}{\tau_i}e_k+T_s\frac{1}{\tau_t}s_k.
\end{alignat}
\end{subequations}
$K_p, \tau_i, \tau_d$ are tuning parameters and $\tau_d/N$ is the time constant for the low pass filter limiting the high-frequency gain of the derivative term \citep{Astrom:Hagglund:1995, Aastrom:Wittenmark:2011}. $\bar{u}$ is the operating point of the MV used in the PID, i.e., $\bar{u}_1 = u_{s,1}$ and $\bar{u}_2 = u_{s,2}$, where $u_s$ is the input operating point used for linearizing the nonlinear model. $\tau_t, s_k$ are integrator windup time-constant and error between commanded and saturated control signal, respectively, and $u_{\min}$ and $u_{\max}$ are the constraints of the control signals. We initialise both PID loops with $I_0 = 0$, $D_{-1} = 0$ and $y_{-1} = y_0$.

\subsubsection{Model-based tuning}
We use the SIMC tuning rules  \citep{Skogestad:2003, Skogestad:Postlethwaite:2005} to tune the SISO PID loops, $(y_1, z_1, u_2)$ and $(y_2, z_2, u_1)$. We apply the transfer functions of the linearized model in \eqref{eq:transfer_functions} to compute the $K_p$, $\tau_i$, and $\tau_d$ in both PID loops as
\begin{equation}
    K_p = \Tilde{K}_p\alpha, \quad\tau_i =  \Tilde{\tau}_i\alpha,\quad \tau_d = \frac{\Tilde{\tau}_d }{\alpha}, 
\end{equation}
where $\alpha = 1+\frac{\Tilde{\tau}_d}{\Tilde{\tau}_i}$ and
\begin{equation}
    \Tilde{K}_p = \frac{\tau_1}{k T_c}, \quad \Tilde{\tau}_i = \min(\tau_1, 4T_c),\quad \Tilde{\tau}_d = \tau_2.
\end{equation}

We choose the tuning parameter $T_c = 50$ and the filter constant $N=5$ for both PID loops in the PID control system. Finally, both PID loops use the integrator windup coefficient $\tau_t = 0.5\tau_i$.

\subsection{Linear model predictive control}
We design an input-bound linear OCP for the LMPC that penalizes tracking errors and input rate of movement, using a discrete-time representation of \eqref{eq:linear_state_space_model} \citep{Sazuan:Joergensen:2018}. We solve the OCP for the LMPC at time $t_k$ with the prediction horizon $T_N$ using the filtered state vector, $\hat{X}_{k|k}$, and disturbance vector, $\hat{D}_{k|k}$, from the CD-KF as initial conditions. We represent the horizon at $t_k$, $[t_k, \, t_k+T_N]$, on an equidistant grid with grid size $N = T_N/T_s$. The disturbances are assumed constant through the horizon. The OCP is
\begin{subequations}
\label{eq:LNMPC_problem}
\begin{alignat}{3}
     &\min_{U,X,Z}\quad && \phi_k = \phi_{Z,k} + \phi_{\Delta U,k} \\
    & s.t.&&
    X_{k} = \hat{X}_{k|k},\\
    & && D_{k+j} = \hat{D}_{k|k}, && j \in \mathcal{N},\\
    & && X_{k+j+1} =
    \Bar{A}X_{k+j}+\Bar{B}U_{k+j}+\Bar{E}D_{k+j}, \quad  && j \in \mathcal{N},\\
    & && Z_{k+j+1} = C_z X_{k+j+1},  && j \in \mathcal{N},\\
    & && U_{\min} \leq U_{k+j} \leq U_{\max}, &&  j \in \mathcal{N},
\end{alignat}
\end{subequations}
where $\mathcal{N} = \{0,1, \ldots, N-1\}$. $\Bar{A}$, $\Bar{B}$, and $\Bar{E}$ are computed from the matrices $A(\theta)$, $B(\theta)$, $E(\theta)$, assuming that the inputs are constant within the sample time, $T_s$. This corresponds to a zero-order hold discretization. The terms in the objective function, $\phi_k$, are
\begin{subequations}
\label{eq:LNMPC_objectives}
\begin{alignat}{3}
     \phi_{Z,k} &= \frac{1}{2}\sum_{j=1}^{N}\left \|   Z_{k+j}-\Bar{Z}_{k+j}\right \|^2_{Q}, \label{subeq:track_LMPC_objective}\\
     \phi_{\Delta U,k} &= \frac{1}{2}\sum_{j=0}^{N-1}\left \|   \Delta U_{k+j}\right \|^2_{S},
\end{alignat}
\end{subequations}
where $\Bar{Z}_{k+j}$ is the setpoint in deviation variables, and $\Delta U_{k+j} = U_{k+j}-U_{k+j-1}$ is the input rate of movement in deviation variables. $Q$ and $S$ are weight matrices. We formulate the OCP as the quadratic program (QP)
\begin{subequations}
\label{eq:sq_problem_LMPC}
\begin{alignat}{3}
        & \min_{\bar{U}_k} \qquad && \phi_k = \frac{1}{2}{\bar U}_k'H {\bar U}_k +g_k' {\bar U}_k + \rho_k\\
    &s.t. &&     \bar{U}_{\min} \leq \bar{U}_k \leq \bar{U}_{\max},
\end{alignat}
\end{subequations}
where $\bar{U}_k = (U_{k+0}, U_{k+1}, U_{k+2}\dots, U_{k+N-1})$ are the decision variables of the QP.
We solve the QP at time $t_k$ to obtain ${\bar U}_k^*$, and implement $u(t_k) = {U}^*_{k+0}+u_s$ to the QTS.
We choose the weight matrices for the LMPC as
\begin{equation}
   Q = \mathrm{diag}\big(\begin{bmatrix}
10 & 10
\end{bmatrix}\big), \quad
        S = \mathrm{diag}\big(\begin{bmatrix}
1 & 1
\end{bmatrix}\big),
\label{eq:MPC_tuning}
\end{equation}
and the number of prediction steps $N = 160$. Given $T_s=5$ s, the prediction horizon of the LMPC is $T_N = 800 \,\mathrm{s} \approx 13 \,\mathrm{min}$.

\subsection{Nonlinear model predictive control}
For the NMPC, we use an input-bound-constrained OCP based on the continuous-time nonlinear model \eqref{eq:model}. We solve the OCP at time $t_k$ with the horizon $T_N$ using the filtered states, $\hat{x}_{k|k}$, and disturbance, $\hat{d}_{k|k}$, from the CD-EKF as initial conditions. We divide the control horizon, $[t_k, \, t_k+T_N]$, into $N = T_N/T_s$ equally spaced subintervals, $[t_{k+j}, \, t_{k+j+1}]$, with  $j \in \mathcal{N} = 0,\dots, N-1$. We assume that the disturbances are constant during the control and prediction horizon. The OCP is
\begin{subequations}
\label{eq:OCP_NMPC}
\begin{alignat}{5}
    &\!\min_{u,x,z}        &\quad& \phi_k = \phi_{z,k} + \phi_{\Delta u,k} &\\
&{s.t.} &      & x(t_k) =\hat{x}_{k|k}, &&\\
&                  &      & d(t) = \hat{d}_{k|k}, &&t\in [t_k, \, t_k+T_N], \\
& & &\dot{x}(t) = f(x(t), u(t), d(t), \theta), \quad &&t\in [t_k, \, t_k+T_N], \\
& & &z(t) = h(x(t), \theta), && &\label{subeq:output_equation} \\
& & &u(t) = u_{k+j}, \quad \quad \quad \quad \, j \in \mathcal{N}, \,\, &&  t\in [t_{k+j}, \, t_{k+j+1}[,\\
& & &u_{\min} \leq u_{k+j} \leq u_{\max}, \quad j\in \mathcal{N},
\end{alignat}
\end{subequations}
and the objective function terms,
\begin{subequations}
\label{eq:objectives_OCP_NMPC}
\begin{alignat}{1}
       \phi_{z,k} &= \frac{1}{2}\sum_{j=1}^{N}\left \|  z_{k+j}-\Bar{z}_{k+j}\right \|_{Q}^2,\\
       \phi_{\Delta u,k} &= \frac{1}{2}\sum_{j=0}^{N-1}\left \|  \Delta u_{k+j}\right \|_{S}^2,
\end{alignat}
\end{subequations}
where $z_{k+j} = z(t_{k+j})$ and $\Bar{z}_{k+j} = \bar{z}(t_{k+j})$.
We transcribe the OCP into a nonlinear program (NLP) using direct multiple-shooting \citep{Bock:Plitt:1984}. The NLP at time $t_k$ is 
\begin{subequations}
\label{eq:NLP_OCP}
\begin{alignat}{5}
      &\min_{w_k} \quad && \phi_k  = \frac{1}{2}\sum_{j=1}^{N}\left \|  z_{k+j}-\Bar{z}_{k+j} \right \|_{Q}^2  + \sum_{j=0}^{N-1}\left \|  \Delta u_{k+j}\right \|_{S}^2 && \\ 
      &s.t. && s_{k+0} = \hat{x}_{k|k},\\
      & && d_{k+j} = \hat{d}_{k|k}, && j \in \mathcal{N},\\
      & && s_{k+j+1} = F_{k+j}(s_{k+j}, u_{k+j}, d_{k+j}, \theta), \quad && j \in \mathcal{N}, \label{subeq:continuity}\\
      & && z_{k+j+1} = h(s_{k+j+1}, \theta), \qquad \quad && j \in \mathcal{N}, \\
      & && u_{\min} \leq u_{k+j} \leq u_{\max}, \qquad \quad && j \in \mathcal{N},
\end{alignat}
\end{subequations}
where $F_{k+j}(s_{k+j}, u_{k+j}, d_{k+j}, \theta)$ represents the state, $x_{k+j+1} = x(t_{k+j+1})$, computed using a fourth-order explicit Runge-Kutta method with 10 fixed integration steps. $w_k = [s_{k+0}; u_{k+0}; s_{k+1}; u_{k+1}; \dots; s_{k+N-1}; u_{k+N-1}; s_{k+N}]$ are the decision variables of the NLP at time $t_k$. We solve the
NLP at time $t_k$ and implement  $u(t_k) = u_{k+0}^*$ on the QTS.

The NMPC uses the open-source software tool CasADi to construct the NLP and we apply IPOPT to solve it \citep{Andersson:etal:2019,Wachter:Biegler:2006}.  At the very first call to OCP at time $t_0 = 0$, we initialize the decision variables with the initial conditions, $ w_0^{[0]} = [x_0; u_0; x_0; u_0; \dots; x_0; u_0; x_0]$. For the subsequent call of the OCP at times $t_{k+1}$ ($k \geq 0$), the previous solution, $w_k^* = [s^*_{k+0}; u^*_{k+0}; s^*_{k+1}; u^*_{k+1}; \ldots; s^*_{k+N-1}; u^*_{k+N-1}; s^*_{k+N}]$, is available, and we apply a shifted expression of the previously converged solution, $w_{k+1}^{[0]} =  [s^{[0]}_{k+1+0}; u^{[0]}_{k+1+0}; s^{[0]}_{k+1+1}; u^{[0]}_{k+1+1}; \ldots; s^{[0]}_{k+1+N-1}; u^{[0]}_{k+1+N-1}; s^{[0]}_{k+1+N}] = [s^{*}_{k+1}; u^*_{k+1}; s^{*}_{k+2}; u^*_{k+2}; \dots; s^*_{k+N}; u^*_{k+N-1}; s^*_{k+N}]$. Note that in such an initialization, $u^{[0]}_{k+N} = u^{[0]}_{k+1+N-1} = u^*_{k+N-1}$ and $s^{[0]}_{k+N+1}=s^{[0]}_{k+1+N} = s^*_{k+N}$, where $u^*_{k+N-1}$ and $s^*_{k+N}$ denote part of the optimal solution $w_k^*$ from the OCP at time $t_k$. The NMPC applies the same weight matrices and number of prediction steps as the LMPC.

NMPC algorithms, that combines an OCP with a CD-EKF based on a stochastic continuous-discrete-time representation of the plant, have successfully been applied to other systems, such as the Van der Pol Oscillator \citep{Brok:2018}, the continuous-stirred tank reactor \citep{Wahlgreen:2020,Jorgensen:etal:2020,Kaysfeld:etal:2023}, spray dryers for milk powder production \citep{Petersen:etal:2017,Petersen:2016}, the U-loop reactor for single-cell protein (SCP) production \citep{Drejer:etal:2017,Ritschel:etal:2019,Ritschel:etal:EKF:2019,Ritschel:2020,Nielsen:2023}, fermentation based bio-manufacturing \citep{Kaysfeld:2023}, chromatography processes \citep{Horsholt:etal:CCTA:2019,Schytt:Jorgensen:NMPC:2024}, isoenergetic-isochoric flash processes \citep{Ritschel:2018,Ritschel:Jorgensen:EKF:2018,Ritschel:Jorgensen:Filters:CCTA:2018,Ritschel:Jorgensen:2019,Ritschel:etal:CompChemEng:2018}, and automated insulin delivery (AID) systems for people with type 1 diabetes \citep{Boiroux:2009,Boiroux:2013,Boiroux:etal:bcell:2010,Boiroux:etal:2010,Boiroux:etal:DYCOPS:2010,Boiroux:etal:ECC:2016,Boiroux:Jorgensen:2017,Boiroux:Jorgensen:NMPC:CDC:2018,Boiroux:Jorgensen:2018,Reenberg:etal:2022,Reenberg:2023,Lindkvist:etal:2023}. Computational and numerical aspects for such problems and their sub-problems have also been addressed \citep{Jorgensen:2004, Jorgensen:2007eccadj,Jorgensen:2007,Kristensen:etal:2004esdirkODEnmpc,Kristensen:etal:2004cace,Kristensen:etal:2005esdirkDAEnmpc,Jorgensen:etal:2018esdirk,Jorgensen:etal:2004,Jorgensen:etal:2012,Capolei:Jorgensen:2012,Frison:2016,Frison:Jorgensen:2013,Frison:Jorgensen:2013cdc,Frison:etal:2014cca,Frison:Sokoler:Jorgensen:2014,Frison:etal:2014,Frison:etal:2016cdc,Frison:etal:2018blasfeo,Frison:Diehl:2020,Edlund:Sokoler:Jorgensen:2009,Sokoler:etal:2013ieeemsc,Sokoler:etal:2013,Sokoler:etal:2016,Christensen:Jorgensen:2024,Jorgensen:etal:2007,Jorgensen:etal:ECC2007ekf,Jorgensen:etal:ACC2007ekf,Boiroux:etal:2019,Binder:etal:2001,Diehl:etal:2009nmpc,Zavala:Biegler:2009,Vershueren:etal:2022,Pulsipher:etal:2022unifying}.

\section{Results and discussion}
\label{sec:results_discussion}

This section presents the experimental results of the three control strategies applied to the physical setup of the QTS. We also provide simulation studies.

\subsection{Experimental results}

We obtain experimental closed-loop data from the physical setup of the QTS. The closed-loop experimental data comprises one experiment for each of the three control strategies, using the same setpoint sequence.

\subsubsection{Software implementation details}
\label{sec:rtapc}

\begin{figure*}[tb]
    \centering
    \includegraphics[width = \textwidth]{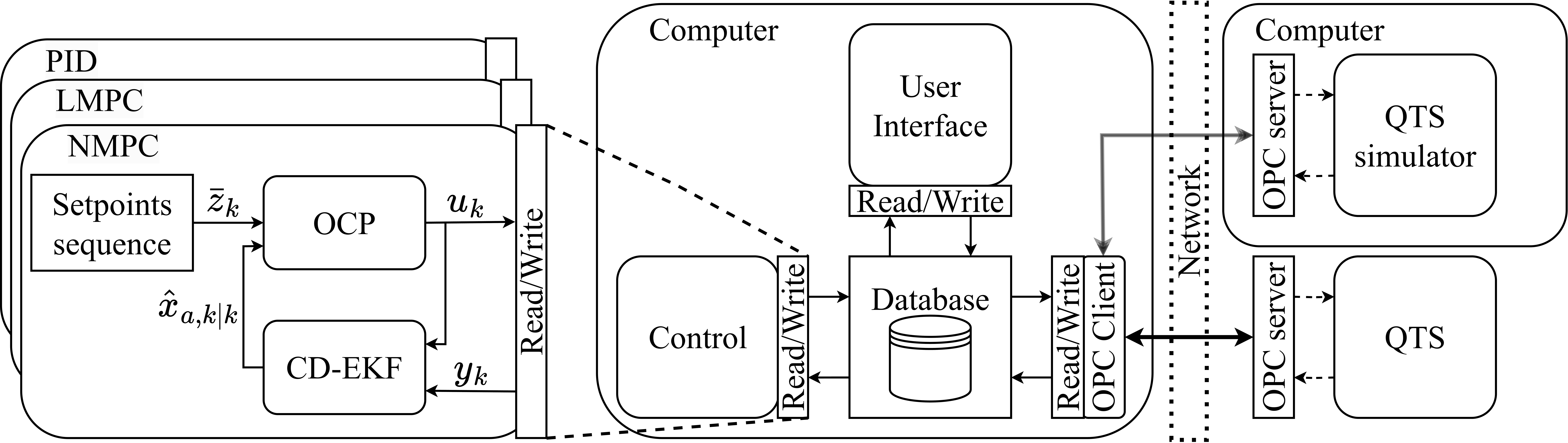}
    \caption{Schematic diagram of the software framework used to implement the three control algorithms (PID, LMPC, NMPC) to the QTS.}
    \label{fig:RTAPC_diagram}
\end{figure*}

We apply the software framework described by \cite{Andersen:etal:ecc:2023} to implement the three control strategies to the physical setup of the QTS. Figure \ref{fig:RTAPC_diagram} illustrates the main components of this framework.
It utilizes an open platform communications unified architecture (OPC UA) connection to communicate with the physical setup of the QTS. It stores measurements and MVs in a PostgreSQL database. The client module periodically transmits data between the database and the QTS using interval timers. 
The control module implements the three control strategies. The computations in these strategies are also performed periodically using interval timers with a period of $T_s$, i.e., the algorithm receives measurement data from the database, computes $u_k$, and writes this value to the database. The user interface module enables operator interaction with the closed-loop system. However, the time-varying setpoints used to compare the three control strategies are pre-programmed to ensure identical experiments. The server-client communication model facilitates real-time simulation experiments by connecting the OPC UA client with a simulator for the QTS. 

\subsubsection{Experimental data}
\label{sec:experiments}

\begin{figure*}[tb]
    \centering
    \centerline{\includegraphics[width=\textwidth]{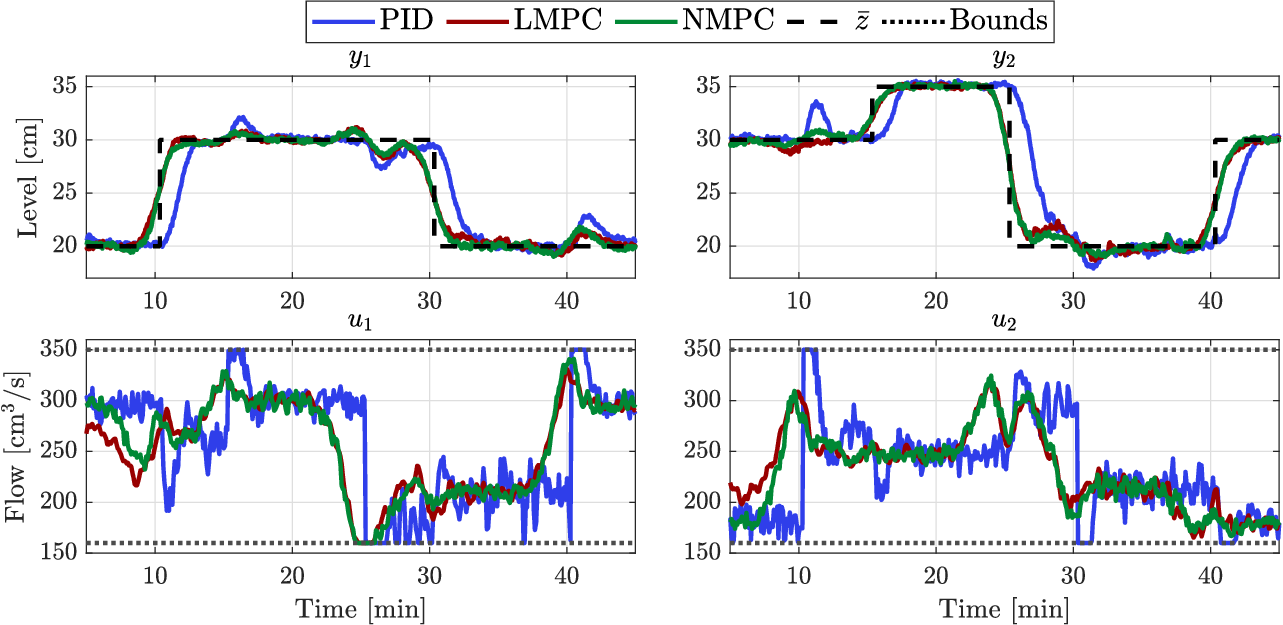}}
    \caption{Data for experiments of decentralized PID, LMPC, and NMPC on the physical QTS.}
    \label{fig:comparison_ML}
\end{figure*}

\begin{figure}[tb]
    \centering
    \centerline{\includegraphics[width=0.5\textwidth ]{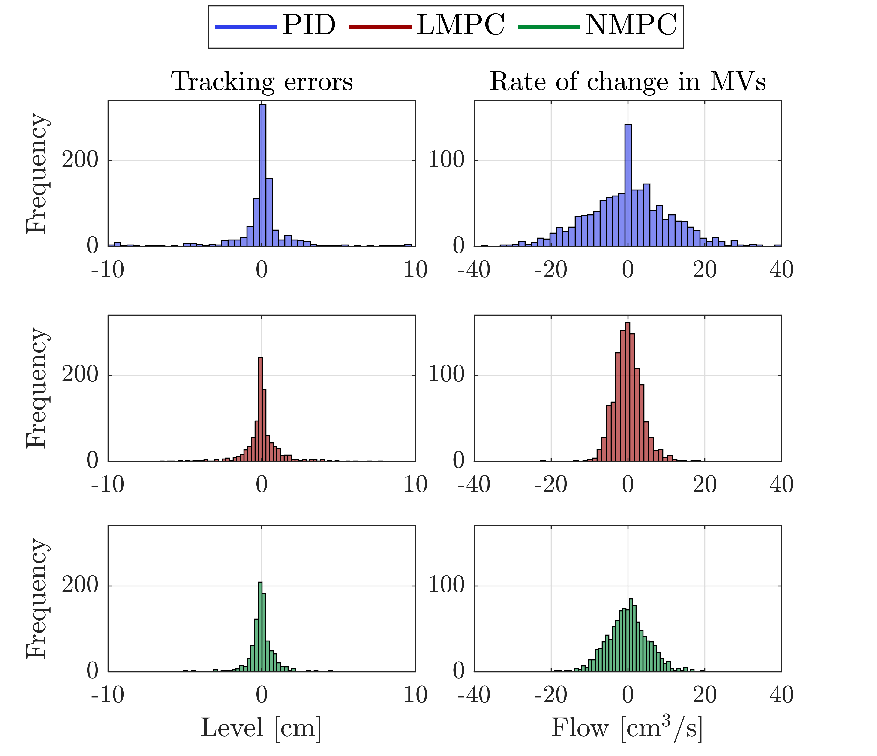}}
    \caption{Histograms of tracking errors and rate of movement in the MVs for the data in Figure \ref{fig:comparison_ML}.}
    \label{fig:histograms}
\end{figure}

We perform one experiment for each of the three control strategies. We apply a sequence of predefined piecewise constant setpoints for both tanks. These setpoints change one by one. The MVs are bounded between $160\: \mathrm{cm}^3/\mathrm{s} \leq u_{i}(t) \leq 350\: \mathrm{cm}^3/\mathrm{s}$ for $i \in \{1, 2\}$. These bounds are motivated by practical considerations related to the physical QTS setup. For the linearization of \eqref{eq:model}, we use the operating point $(x_s, u_s, d_s)$ by choosing $u_s~=~\begin{bmatrix}
300\:\mathrm{cm}^3/\mathrm{s}; & 300\:\mathrm{cm}^3/\mathrm{s}
\end{bmatrix}$, $d_s$ as the zero vector, 
and obtain $x_s$ by solving $0~=~f(x_s, u_s, d_s, \theta_{ML}^*)$. 

We initialize the experiments by letting the controller reach the initial phase of the setpoints and waiting for a steady-state in the QTS. Subsequently, the setpoint sequence is started. Figure \ref{fig:comparison_ML} presents the data from the experiments and
Figure \ref{fig:histograms} shows histograms of the tracking errors and rate of movement in the MVs.


\subsection{Simulation studies}
\label{sec:simulations}

This section presents four simulation case studies of the three control strategies. The two first simulation studies apply the same sequence of setpoints as in Figure \ref{fig:comparison_ML}. In these simulations, we compare the influence of providing future setpoints to the MPCs by simulating them with and without future setpoint information. For the situation without future setpoint information, the setpoints are assumed constant over the prediction horizon for both MPCs. In the third simulation study, we simulate the closed-loop systems affected by large deterministic disturbances in tank 3 and tank 4. For these simulations, we keep the setpoints constant. In the fourth simulation study, we apply a large diffusion term to excite the closed-loop systems with stochastic disturbances. We also keep the setpoints constant in these simulations. We apply the nominal parameters in Table \ref{tab:all_parameters} for all four simulation studies.

\subsubsection{Tuning of controllers}
For the simulation case studies, we use the same tuning for the decentralized PID control system and the same weight matrices and prediction horizon for the MPCs. The CD-KF and CD-EKF, however, are tuned with $\sigma_{a,ii}(\theta) = 1.0$ for $i=1, \dots, 8$, and we apply the measurement noise covariance $R = \text{diag}([0.02, 0.02, 0.02, 0.02])$. We also use these covariances for simulating \eqref{eq:model}. As a result, the CD-KF and CD-EKF have perfect knowledge about the systems. For the fourth simulation study, we choose
\begin{equation}
    \sigma(\theta) = \text{diag}\big(\begin{bmatrix}
 20 & 20 & 20 & 20 
\end{bmatrix}\big)
\end{equation} 
and apply the diffusion coefficients 
\begin{equation}
    \sigma_a(\theta) = \text{diag}\big(\begin{bmatrix}
0 & 0 & 0 & 0 & 20 & 20 & 20 & 20 
\end{bmatrix}\big)
\end{equation}
in the CD-KF and CD-EKF.

\begin{figure}[tb]
    \centering
    \centerline{\includegraphics[width=0.5\textwidth]{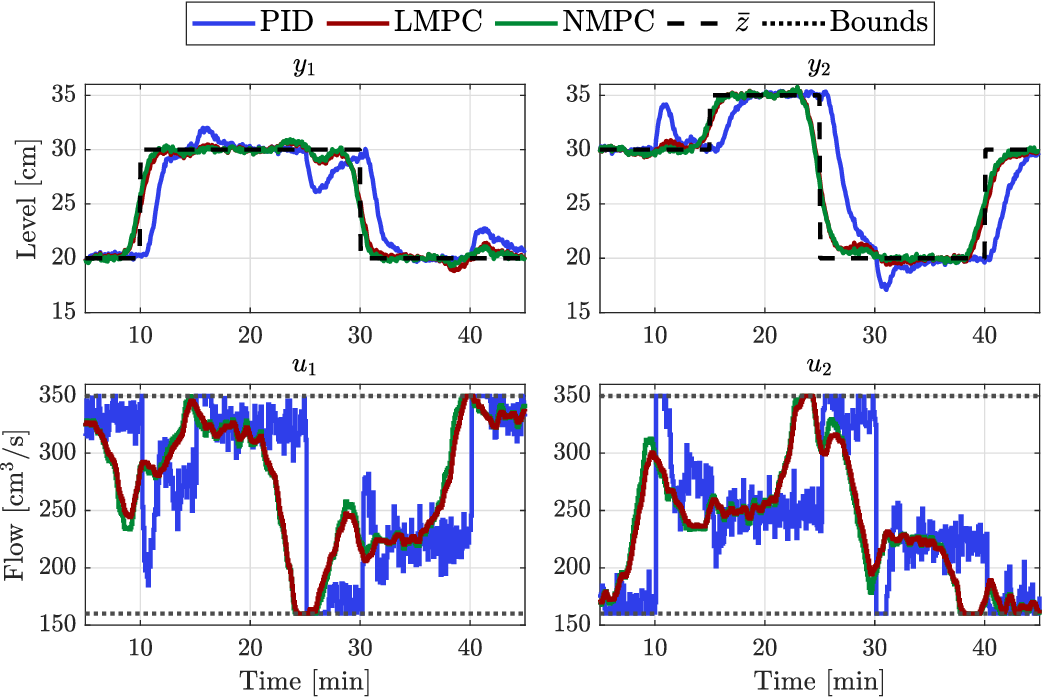}}
    \caption{Simulation 1: Tracking predefined time-varying setpoints. The MPCs receive information on future setpoint changes.}
\label{fig:sim_anticipatory_vs_reactive_v2_a}
\end{figure}

\begin{figure}[tb]
    \centering
    \centerline{\includegraphics[width=0.5\textwidth]{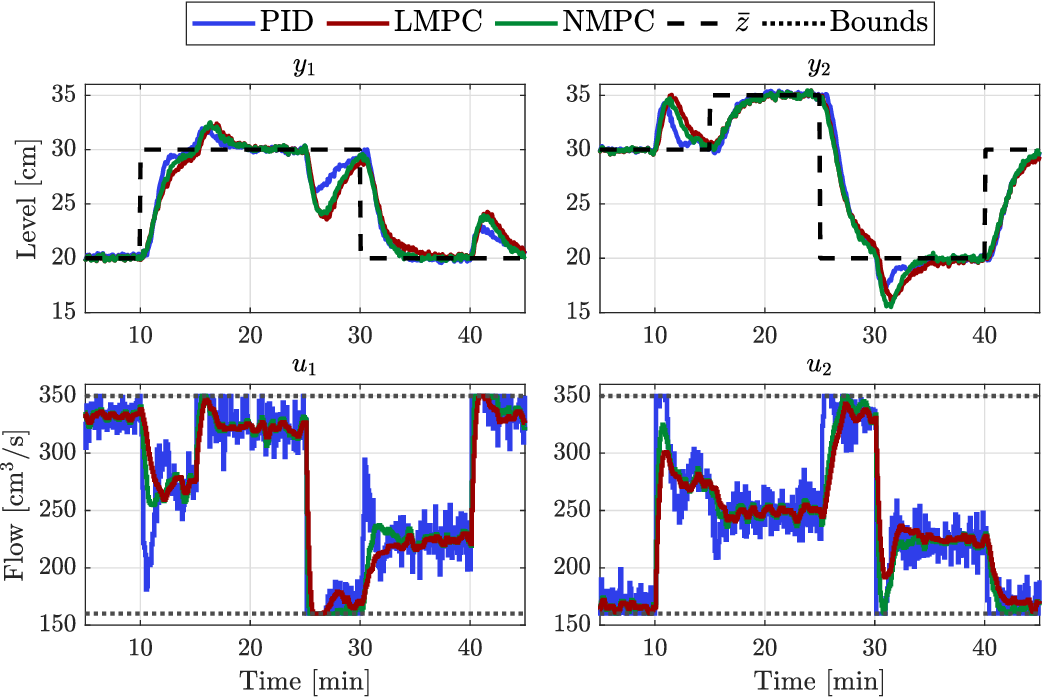}}
    \caption{Simulation 2: Tracking predefined time-varying setpoints. The MPCs receive only information about current setpoint values.}
\label{fig:sim_anticipatory_vs_reactive_v2_r}
\end{figure}

\begin{figure}[tb]
    \centering
    \centerline{\includegraphics[width=0.5\textwidth]{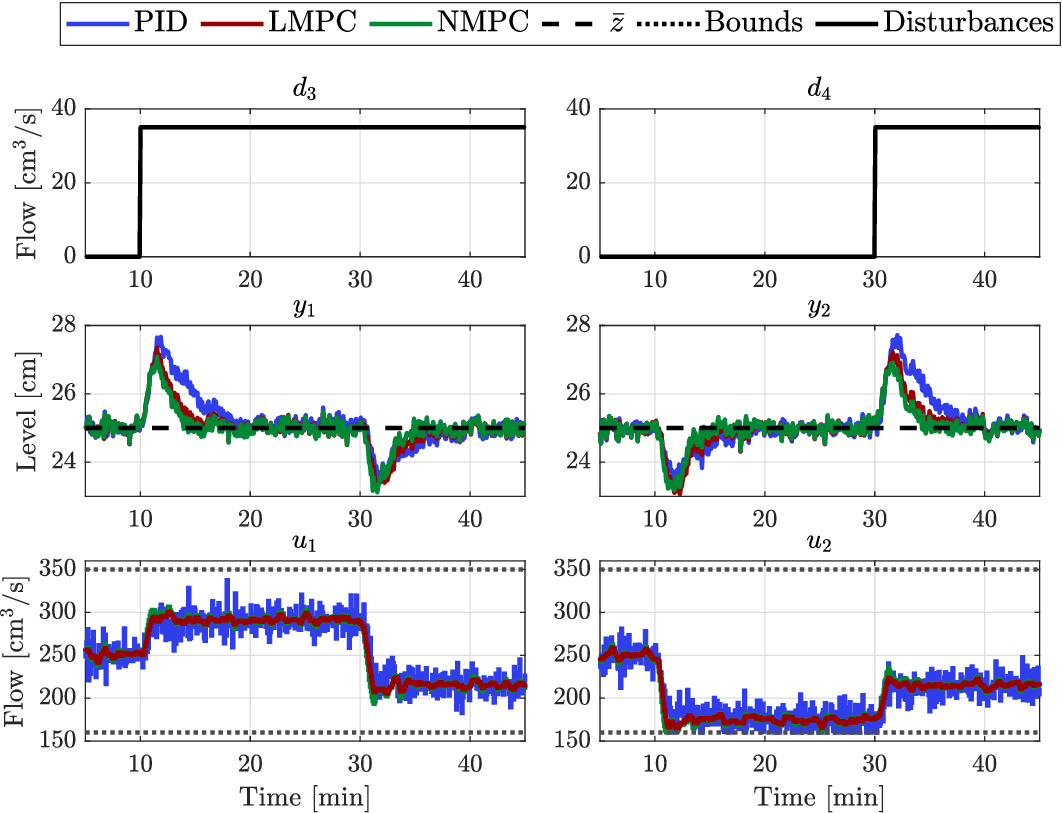}}
    \caption{Simulation 3: Constant setpoints with large deterministic disturbance step-changes in Tank 3 and Tank 4.}
\label{fig:sim_anticipatory_vs_reactive_v3_a}
\end{figure}

\begin{figure}[tb]
    \centering
    \centerline{\includegraphics[width=0.5\textwidth]{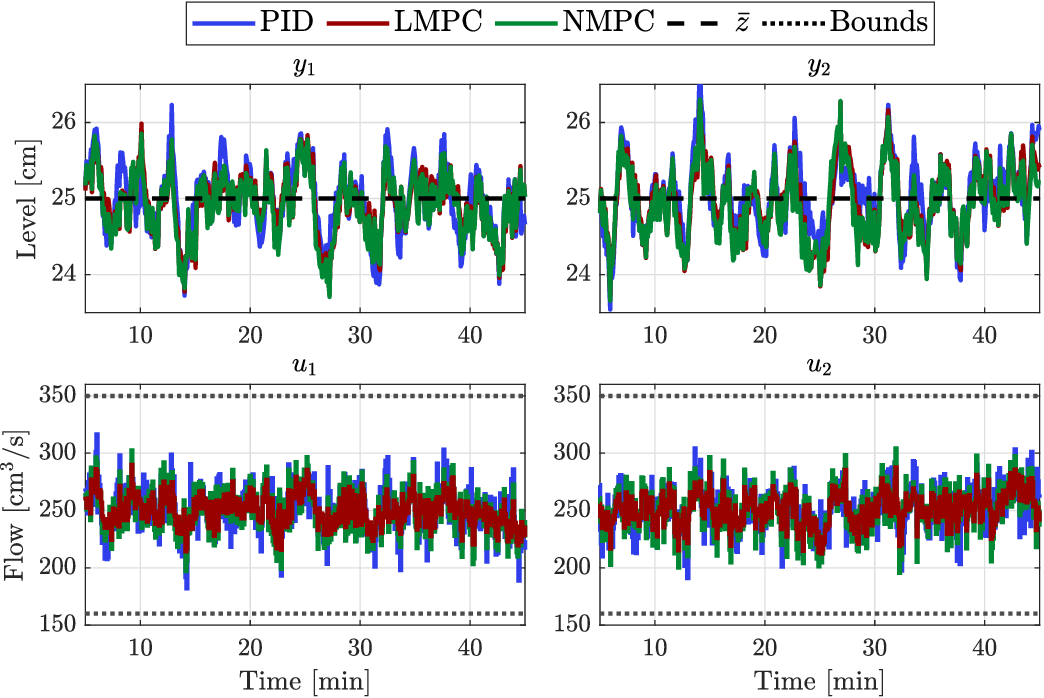}}
    \caption{Simulation 4: Constant setpoints and no deterministic disturbance step changes. The QTS is simulated using a large diffusion term.}
\label{fig:sim_anticipatory_vs_reactive_v4_a}
\end{figure}

Figures \ref{fig:sim_anticipatory_vs_reactive_v2_a}-\ref{fig:sim_anticipatory_vs_reactive_v2_r} present the simulations for setpoint tracking with and without future setpoint information in the MPCs. In Figure \ref{fig:sim_anticipatory_vs_reactive_v3_a}, we demonstrate disturbance rejection of large deterministic disturbances and for constant setpoints. Finally, Figure \ref{fig:sim_anticipatory_vs_reactive_v4_a} shows the controller performances for constant setpoints and no deterministic disturbance step-changes. In these results, we simulate a large diffusion term to excite the system with stochastic disturbances.

\subsubsection{Evaluation of controller performances}

We measure the performance of the controllers by applying different norms to the observed tracking error and rate of movement in the MVs. The performance measures are the normalized integral squared error (NISE), the normalized integral absolute error (NIAE), and the normalized integral squared rate of movement in the MVs (NIS$\Delta$U). These performance measures are defined as
\begin{subequations}
\label{eq:performance_metrics}
\begin{alignat}{1}
         \mathrm{NISE} &= \frac{1}{N}\sum_{k = 1}^N\left \|   \bar{e}_k\right \|^2_2, \\
         \mathrm{NIAE} &= \frac{1}{N}\sum_{k = 1}^N\left \|   \bar{e}_k\right \|_1,\\
       \mathrm{NIS\Delta U} &= \frac{1}{M-1}\sum_{k = 1}^{M-1}\left \|   \Delta u_k\right \|^2_2,
\end{alignat}
\end{subequations}
where $N$ and $M$ are the number of data points of the water level and the MVs, respectively, and $\bar{e}_k$ is the observed tracking error
\begin{equation}
    \bar{e}_k = \bar{z}_k - \begin{bmatrix} y_{1,k}; & y_{2,k}\end{bmatrix}.
\end{equation}
We apply these performance measures to the experimental data and the data of the four simulation studies. Tables \ref{table:setpoint_performance_metrics}-\ref{table:dist_rej_performance_metrics} present the results.

\begin{table*}[tb]
\centering
\caption{Performance measures of decentralized PID, LMPC, and NMPC for tracking predefined time-varying setpoints. Experimental data and Simulation 1 apply MPCs that receive information on future setpoint changes. In Simulation 2, only current setpoint information is provided for the MPCs. Bold numbers indicate the lowest value in their respective categories.}
\label{table:setpoint_performance_metrics}
\begin{tabular}{llllllllllll}
\hline
          & \multicolumn{3}{c}{Experimental (Figure \ref{fig:comparison_ML})} &  & \multicolumn{3}{c}{Simulation 1 (Figure \ref{fig:sim_anticipatory_vs_reactive_v2_a})} &  & \multicolumn{3}{c}{Simulation 2 (Figure \ref{fig:sim_anticipatory_vs_reactive_v2_r})} \\ \cline{2-4} \cline{6-8} \cline{10-12} 
          & NISE   & NIAE   & NIS$\Delta$U   &  & NISE   & NIAE   & NIS$\Delta$U   &  & NISE   & NIAE   & NIS$\Delta$U   \\
PID       &   9.063 & 1.459 & 249.079 &  &  14.214      & 3.480       &  784.106      &  &\textbf{14.214}      & \textbf{3.480}       &  784.106        \\
LMPC      &   1.637 & 0.728 & \textbf{12.089} &  &   2.232      &   1.391     & \textbf{11.003}      &  &   18.688     &   4.111     &   \textbf{38.210}     \\
NMPC      &  \textbf{1.423} & \textbf{0.647} & 28.674 &  &   \textbf{1.941}    &    \textbf{1.222}    &   18.262   &   & 17.194 &       3.775  &   81.787           \\ \hline
\end{tabular}
\end{table*}

\begin{table*}[tb]
\centering
\caption{Decentralized PID, LMPC, and NMPC performance measures for disturbance rejection and constant setpoints. Simulation 3 involves deterministic disturbance rejection. Simulation 4 is the simulation study using a large diffusion term. Bold numbers indicate the lowest value in their respective categories.}
\label{table:dist_rej_performance_metrics}
\begin{tabular}{llllllll}
\hline
          & \multicolumn{3}{c}{Simulation 3 (Figure \ref{fig:sim_anticipatory_vs_reactive_v3_a})} &  & \multicolumn{3}{c}{Simulation 4 (Figure \ref{fig:sim_anticipatory_vs_reactive_v4_a})} \\ \cline{2-4} \cline{6-8} 
          & NISE   & NIAE   & NIS$\Delta$U   &  & NISE   & NIAE   & NIS$\Delta$U   \\
PID       &   1.046     &    0.896    &     740.750    &  &  0.294     &  0.603     &  962.480     \\
LMPC      &    0.670    &    0.673    &    \textbf{7.206}      &  &   \textbf{0.264}     &     \textbf{0.595}
     &    \textbf{287.044}    \\
NMPC      &   \textbf{0.555}     &      \textbf{0.593}  &    15.659    &  &      0.280  &     0.606    &   834.926    \\ \hline
\end{tabular}
\end{table*}

\subsection{Discussion}
This section discusses the experimental and simulated results and compares these in terms of setpoint tracking, rate of movement in the MVs, and disturbance rejection capabilities. 

\subsubsection{Influence of future setpoint information in the MPCs}

The experimental results in Figure \ref{fig:comparison_ML} show improved tracking capabilities of the LMPC and the NMPC compared with the decentralized PID control system, for predefined time-varying setpoints. Figure \ref{fig:histograms} also shows that the outliers in the tracking errors are removed when using MPCs instead of the decentralized PID control system. The simulations in Figures \ref{fig:sim_anticipatory_vs_reactive_v2_a}-\ref{fig:sim_anticipatory_vs_reactive_v2_r} show that the MPCs without future setpoint information perform worse than the MPCs with future setpoint information. If the future setpoint information is removed from the MPCs, the decentralized PID tracking errors appear to be even better than the MPCs. This is demonstrated in Table \ref{table:setpoint_performance_metrics}, where NISE and NIAE are smaller for the decentralized PID control system than for the LMPC and the NMPC.

\subsubsection{Rate of movement in the MVs}

Besides improved setpoint tracking, the experimental results also show a significant reduction in the rate of movement of the MVs of the MPCs compared with the decentralized PID. This is clearly illustrated in Figure \ref{fig:histograms} that shows the squeezed distributions of the rate of movement in the MVs when using MPCs instead of decentralized PID. This is assumed to be due to the proportional part of the PIDs, which reacts very rapidly to setpoint changes by switching between the input bounds. This behavior is also the main contributor to the PID controller's fast response, which leads to good tracking.

\subsubsection{Disturbance rejection}

Table \ref{table:dist_rej_performance_metrics} summarizes the simulation results for two different disturbance scenarios. For large deterministic step changes in the disturbances  (Simulation 3), the MPCs perform slightly better than the decentralized PID in terms of NISE and NIAE. 
In Simulation 4, the decentralized PID and the MPCs perform almost identically regarding setpoint tracking. The only notable performance difference between the three controllers is in the NIS$\Delta$U where the LMPC performs better than the decentralized PID and the NMPC.

\subsubsection{Comparison between the LMPC and the NMPC}

We see no notable differences in performance between the LMPC and the NMPC for the experimental results.  Futhermore, the experimental results are consistent with the four simulation studies. Here, the NMPC performs slightly better in terms of NISE and NIAE. However, the NIS$\Delta$U is significantly higher for the NMPC when compared with the LMPC. 

We assume that the main reason for the negligible difference between the LMPC and the NMPC is the degree of nonlinearity of the QTS. The QTS is a nonlinear system, but for the range of operation, the linearized model approximates the nonlinear system well. The differences in NISE and NIAE for LMPC and NMPC are slightly bigger in Simulations 1 and 2 compared to Simulations 3 and 4, as the range of operations in Simulations 1 and 2 are larger compared to the others. As a result, the linearized model is a better approximation of the nonlinear system in Simulations 3 and 4.

\subsubsection{Limitations of the comparison between PID and MPC}
The experimental and simulated results were performed using only one set of tuning parameters for each controller. Obviously, the results may have been different if alternative tunings for the LMPC and NMPC (weight matrices, prediction horizon) as well as the PIDs (closed-loop time constants, $T_c$) were chosen. The tuning parameters for the PID type controller were based on model-based tuning rules i.e. SIMC, while the tuning of the MPC were based on trial-and-error adjustment of the weights. In this way, each controller consists of its algorithm as well as its method for tuning. Often the tuning of MPCs is less addressed than detailed algorithmic description \cite[chap. 7]{Maciejowski:2002}. 

While significant effort was made to obtain a tuning that makes the comparison fair by providing the best possible practical tuning for each controller, the comparison could be augmented and improved by using model- and optimization-based tuning for the PID controller as well as the MPCs (LMPC and NMPC). Tuning the decentralized PID controller using an optimization-based approach would be another possibility to investigate the best achievable performance of a PID controller. Especially if the objective functions chosen are identical to those of the MPCs. The performance measure \eqref{eq:performance_metrics} favors the MPCs because it is based on the same objectives as in the MPCs. In an uncertainty quantification study, the PID controller as well as the MPCs could be tuned by optimization-based approaches to have the best possible behavior of all controllers that were to be compared \citep{Olesen:etal:2013siso,Olesen:etal:2013mimo}. Such systematic computational brute-force tuning approaches may be a subject for further studies to compare the inherent performance of PID and MPC algorithms \citep{Wahlgreen:etal:2021,Kaysfeld:etal:2023,Kaysfeld:2023}.  

\section{Conclusions}
\label{sec:conclusions}

This paper compared the performance of a decentralized proportional-integral-derivative (PID) control system, a linear model predictive control (LMPC) strategy, and a nonlinear model predictive control (NMPC) algorithm applied to a quadruple tank system (QTS). Both experimental and simulated results were considered. We modeled the QTS as a nonlinear stochastic continuous-discrete-time system, and we identified the parameters using a maximum-likelihood (ML) prediction-error-method (PEM). The NMPC combined a nonlinear optimal control problem (OCP) with a continuous-discrete extended Kalman filter (CD-EKF) based on the stochastic nonlinear model. The LMPC combined a linear OCP with a continuous-discrete Kalman filter (CD-KF). The linear OCP and the CD-KF applied a linearized version of the nonlinear model. The decentralized PID controller was tuned using the simple internal model control (SIMC) rules. These rules required transfer functions of the QTS, and we obtained these from the linearized model.

We systematically tested the three control systems using predefined time-varying setpoints for the water levels in the two bottom tanks of the QTS. We did this both for a physical setup of the QTS and in simulations. For the physical setup, we provided the future setpoints to the MPCs. In the simulation studies, we repeated the experimental results and compared them with simulations of the MPCs without providing information on future setpoints. We also simulated the controllers' disturbance rejection capabilities.

We evaluated the performance of the control systems in terms of norms of the tracking errors and the rates of movement in the manipulated variables (MVs). The experimental results showed that the LMPC and the NMPC performed better than the decentralized PID control system, while the LMPC and the NMPC performed similarly in terms of norms of the setpoint deviations. However, the simulation studies showed, that the main advantage of the MPCs was their ability to use predefined setpoints. By providing only the current setpoint to the MPCs, the decentralized PID control system achieved better tracking performance. For the rejection of large deterministic disturbances, the MPCs showed slightly improved performance compared to the decentralized PID. However, for rejection of stochastic disturbances, no notable difference between the decentralized PID and the MPCs was observed. Nevertheless, the MPCs still achieved a much lower input rate of movement in all four simulations and experimental results.

We discussed the limitations of the comparisons. The PID tuning was based on the SIMC tuning rules, while the MPCs were tuned by trial-and-error. Computationally intensive optimization based tuning using uncertainty quantification for the MPCs as well as the PID controller would represent the best possible performance for each control algorithm given an identified model. We did not provide such results in the paper, as the computational tools were in development, when we conducted the experimental studies reported in this paper \citep{Wahlgreen:etal:2021,Kaysfeld:etal:2023,Kaysfeld:2023}. However, systematic optimization based tuning of the PID controller and the MPCs for the QTS using uncertainty quantification represents an interesting future study. We also note that the objective functions of the MPCs in this paper are formulated as discrete-time sums, as is common in the linear MPC literature. It would be interesting to implement the LMPC and NMPC using integrals (continuous-time with piecewise constant inputs) as this may provide more intuitive tuning procedures \citep{Zhang:2024,Zhang:etal:ADCHEM2023:2024,Zhang:etal:ECC2024:2024,Zhang:etal:DYCOPS2025:2025,Zhang:etal:ECC2025:2025}.


\section*{Acknowledgements}
This project has been partly funded by the MissionGreenFuels project
DynFlex under The Innovation Fund Denmark no. 1150-00001B, the INNO-CCUS project NewCement under The Innovation Fund Denmark no. 1150-00001B, the Danish Energy Technology and Demonstration Program (EUDP) under the Danish Energy Agency in the EcoClay project no. 64021-7009, and
the IMI2/EU/EFPIA Joint Undertaking Inno4Vac no. 101007799. This communication reflects the authors’ views and that neither IMI nor the European Union, EFPIA, or
any Associated Partners are responsible for any use that may be made of the information contained therein.



\bibliographystyle{elsarticle-harv}   
\bibliography{references_qts_2024}

\end{document}